\documentclass[12pt,fullpage]{article}

\usepackage{amsfonts}
\usepackage{xr}
\usepackage{graphicx}
\usepackage{amsmath}
\usepackage{epsfig}

\def\Z{\mathbb{Z}}
\def\R{\mathbb{R}}
\def\N{\mathbb{N}}

\def\epsilon{\varepsilon}
\def\trait (#1) (#2) (#3){\vrule width #1pt height #2pt depth #3pt}
\def\fin{\hfill\trait (0.1) (5) (0) \trait (5) (0.1) (0) \kern-5pt 
\trait (5) (5) (-4.9) \trait (0.1) (5) (0)}

\usepackage{color}

\definecolor{gr}{rgb}   {0.,   0.8,   0. } 
\definecolor{bl}{rgb}   {0.,   0.5,   1. } 
\definecolor{mg}{rgb}   {0.7,  0.,    0.7}

\newcommand{\SE}{\setcounter{equation}{0} \section}
\newcommand{\be}{\begin{equation}}
\newcommand{\ee}{\end{equation}}
\newcommand{\baa}{\begin{array}}
\newcommand{\eaa}{\end{array}}
\newcommand{\ba}{\begin{eqnarray}}
\newcommand{\ea}{\end{eqnarray}}

\newtheorem{theo}{\bf Theorem}[section]
\newtheorem{lem}[theo]{\bf Lemma}
\newtheorem{pro}[theo]{\bf Proposition}
\newtheorem{cor}[theo]{\bf Corollary}

\newtheorem{rem}[theo]{\bf Remark}

\DeclareMathOperator{\supp}{supp}

\begin{document}
\author{\begin{tabular}{ccc}
Nadine Badr\footnote{Universit\'e de Paris-Sud, Orsay et CNRS UMR 8628, 91405 Orsay Cedex (France) Email: nadine.badr@math.u-psud.fr} &
Emmanuel Russ \footnote{Universit\'e Paul C\'ezanne, Facult\'e des Sciences et Techniques de
Saint-J\'er\^ome, Avenue Escadrille Normandie-Ni\'emen,
13397 MARSEILLE Cedex 20, FRANCE, and LATP, CNRS, UMR 6632. E-mail:
emmanuel.russ@univ-cezanne.fr}\\ \\
Universit\'e de Paris-Sud & Universit\'e Paul C\'ezanne 
\end{tabular}}

\title{Interpolation of Sobolev spaces, Littlewood-Paley inequalities and Riesz transforms on graphs}
\maketitle

\noindent{\small{{\bf Abstract.} }} Let $\Gamma$ be a graph endowed 
with a reversible Markov kernel $p$, and $P$ the associated operator, defined by $Pf(x)=\sum_y p(x,y)f(y)$. Denote by $\nabla$ the discrete gradient. We give necessary and/or sufficient conditions on 
$\Gamma$ in order to compare $\left\Vert \nabla f\right\Vert_{p}$ 
and $\left\Vert (I-P)^{1/2}f\right\Vert_{p}$ uniformly in $f$ for 
$1<p<+\infty$. These conditions are different for $p<2$ and $p>2$. 
The proofs rely on recent techniques developed to handle operators 
beyond the class of Calder\'on-Zygmund operators. For our purpose, we also prove Littlewood-Paley inequalities and interpolation results for Sobolev spaces in this context, which are of independent interest.

\noindent{\small{{\bf AMS numbers 2000: }}} Primary: 60J10. Secondary:  42B20, 42B25.


\tableofcontents
%%%%%%%%%%%%%%%%%%%%%%%%%%%%%%%%%%%%%%%%%%%%%%%%%%%%%%%%%%%%%%%%%%%%%%%%%%%%%%%

%%%%%%%%%%%%%%%%%%%%%%%%%%%%%%%%%%%%%%%%%%%%%%%%%%%%%%%%%%%%%%%%%%%%%%%%%%%%%%%
\SE{Introduction and results}
It is well-known that, if $n\geq 1$, $\left\Vert \nabla 
f\right\Vert_{L^p(\R^n)}$ and $\left\Vert 
(-\Delta)^{1/2}f\right\Vert_{L^p(\R^n)}$ are 
comparable uniformly in $f$ for all $1<p<+\infty$. This fact means that 
the classical Sobolev space $W^{1,p}(\R^n)$ defined by means of the 
gradient coincides with the Sobolev space defined through the Laplace operator. 
This is interesting in particular because $\nabla$ is a local operator, 
while $(-\Delta)^{1/2}$ is not.

Generalizations of this result to geometric contexts can be given. 
On a Riemannian manifold $M$, it was asked by Strichartz in 
\cite{stri} whether, if $1<p<+\infty$, there exists $C_{p}>0$ such 
that, for all function $f\in C^{\infty}_{0}(M)$, 
\begin{equation} \label{equivalence}
C_{p}^{-1}\left\Vert \Delta^{1/2}f\right\Vert_{p} \leq \left\Vert \left\vert df\right\vert\right\Vert_{p}\leq 
C_{p}\left\Vert \Delta^{1/2}f\right\Vert_{p},
\end{equation}
where $\Delta$ stands for the Laplace-Beltrami operator on $M$ and $d$ for the exterior differential. Under suitable assumptions on $M$, which can be formulated, for instance, in terms 
of the volume growth of balls in $M$, uniform $L^{2}$ Poincar\'e inequalities 
on balls of $M$, estimates on the heat semigroup ({\it i.e.} the 
semigroup generated by $\Delta$) or the Ricci curvature, each of the two inequalities 
contained in (\ref{equivalence}) holds for a range of $p$'s (which is, 
in general, different for the two inequalities). The second inequality in (\ref{equivalence}) means that the Riesz transform $d\Delta^{-1/2}$ is $L^p$-bounded. We refer 
to (\cite{ac,acdh,bakry,cd}) and the references therein. \par
In the present paper, we consider a graph equipped 
with a discrete gradient and a discrete Laplacian and investigate the corresponding counterpart 
of (\ref{equivalence}). To that purpose, we prove, among other things, an interpolation result for Sobolev spaces defined via the differential, similar to those already considered in \cite{ostr}, as well as $L^p$ bounds for Littlewood-Paley functionals.  
% $\left\Vert \left\vert\nabla f\right\vert \right\Vert_{p}$ 
% and $\left\Vert \Delta^{1/2}f\right\Vert_{p}$ are comparable uniformly 
% in $f$ for $1<p<+\infty$ where $\Delta$ is the Laplace-Beltrami 
% operator. In \cite{alexo}, Alexopoulos proves that, on a Lie group $G$ equipped with a system of left-invariant vector fields 
% $(X_1,...,X_k)$ satisfying the H\"ormander condition, if $\Delta=-\sum_{1\leq i\leq k} X_i^2$, then 
% $\left\Vert \Delta^{1/2}f\right\Vert_p\sim \sum_{1\leq i\leq k} \left\Vert X_if\right\Vert_p$ for all $1<p<+\infty$. 
% It is shown in \cite{bakry} that, on a Riemannian manifold with non-negative Ricci curvature, 
% $\left\Vert \Delta^{1/2}f\right\Vert_p\sim \left\Vert \nabla 
% f\right\Vert_p$ for all $1<p<+\infty$, where $\Delta$ is the Laplace-Beltrami operator. In \cite{acdh}, a sharp sufficient condition on a Riemannian manifold $M$ is given to ensure that $\left\Vert \nabla f\right\Vert_p$ and $\left\Vert \Delta^{1/2}f\right\Vert_p$ are comparable for all $1<p<+\infty$. This last result encompasses Alexopoulos' and Bakry's statements. In the present paper, we give a discrete version of the main theorem in \cite{acdh}. In other words, if $\Gamma$ is a graph, $\Delta$ a discrete ``Laplace operator'' on $\Gamma$ and $\nabla$ a discrete ``length of gradient'' on $\Gamma$, is it true that $\Delta^{1/2}f$ and $\nabla f$ have comparable $L^p$-norms for all $1<p<+\infty$ ? 
\subsection{Presentation of the discrete framework}
Let us give precise definitions of our framework. The following presentation is partly borrowed from \cite{delmo}. Let $\Gamma$ be an infinite set and 
$\mu_{xy}=\mu_{yx}\geq 0$ a symmetric weight on $\Gamma\times \Gamma$. We call $(\Gamma,\mu)$ a weighted graph. In the sequel, we write most of the time $\Gamma$ instead of $(\Gamma,\mu)$, somewhat abusively. If $x,y\in \Gamma$, say that $x\sim y$ if and only 
if $\mu_{xy}>0$. Denote by $E$ the set of edges in $\Gamma$, {\it i.e.}
\[
E=\left\{(x,y)\in \Gamma\times \Gamma;\ x\sim y\right\},
\]
and notice that, due to the symmetry of $\mu$, $(x,y)\in E$ if and only if $(y,x)\in E$.\par
For $x,y\in \Gamma$, a path joining $x$ to $y$ is a finite sequence of edges $x_0=x,...,x_N=y$ such that, for all $0\leq i\leq N-1$, 
$x_i\sim x_{i+1}$. By definition, the length of such a path is $N$. Assume that $\Gamma$ is connected, which means that, for all $x,y\in \Gamma$, there exists a path joining 
$x$ to $y$. For all $x,y\in \Gamma$, the distance between $x$ and $y$, 
denoted by $d(x,y)$, is the shortest length of a path joining $x$ and $y$. For all $x\in \Gamma$ and all $r\geq 0$, let 
$B(x,r)=\left\{y\in \Gamma,\ d(y,x)\leq r\right\}$. In the sequel, 
we always assume that $\Gamma$ is locally uniformly finite, which means 
that there exists $N\in \N^{\ast}$ such that, for all $x\in \Gamma$, $\sharp 
B(x,1)\leq N$(here and after, $\sharp A$ denotes the cardinal of any subset $A$ of $\Gamma$). If $B=B(x,r)$ is a ball, set $\alpha B=B(x,\alpha r)$ for all $\alpha>0$, and write $C_1(B)=4B$ and $C_j(B)=2^{j+1}B\setminus 2^{j}B$ for all integer $j\geq 2$.\par
\noindent For any subset $A\subset \Gamma$, set
\[
\partial A=\left\{x\in A;\ \exists y\sim x,\ y\notin A\right\}.
\]
For all $x\in \Gamma$, set $m(x)=\sum\limits_{y\sim x}\mu_{xy}$. We always 
assume in the sequel that $m(x)>0$ for all $x\in \Gamma$. If 
$A\subset\Gamma$, define $m(A)=\sum\limits_{x\in A} m(x)$. For all $x\in \Gamma$ and $r>0$, write $V(x,r)$ instead of $m(B(x,r))$ and, if $B$ is a ball, $m(B)$ will be denoted by $V(B)$.\par

For all $1\leq p<+\infty$, say that a function $f$ on $\Gamma$ belongs to $L^p(\Gamma,m)$ (or $L^p(\Gamma)$) if
\[
\left\Vert f\right\Vert_p:=\left(\sum\limits_{x\in \Gamma}\left\vert f(x)\right\vert^pm(x)\right)^{1/p}<+\infty.
\]
Say that $f\in L^{\infty}(\Gamma,m)$ (or $L^{\infty}(\Gamma)$) if
\[
\left\Vert f\right\Vert_{\infty}:=\sup\limits_{x\in \Gamma} \left\vert f(x)\right\vert<+\infty.
\]
Define $p(x,y)=\mu_{xy}/m(x)$ for all $x,y\in \Gamma$. Observe that $p(x,y)=0$ if $d(x,y)\geq 2$. Set also
\[
p_0(x,y)=\delta(x,y)
\]
and, for all $k\in \N$ and all $x,y\in \Gamma$,
\[
p_{k+1}(x,y)=\sum\limits_{z\in \Gamma}p(x,z)p_k(z,y).
\]
The $p_k$'s are called the iterates of $p$. Notice that, for all $x\in \Gamma$, there are at most $N$ non-zero terms in 
this sum. Observe also that, for all $x\in \Gamma$,
\begin{equation} \label{stochastic}
\sum\limits_{y\in\Gamma}p(x,y)=1
\end{equation}
and, for all $x,y\in \Gamma$,
\begin{equation} \label{reversibility}
p(x,y)m(x)=p(y,x)m(y).
\end{equation}
For all function $f$ on $\Gamma$ and all $x\in \Gamma$, define
\[
Pf(x)=\sum\limits_{y\in \Gamma}p(x,y)f(y)
\]
(again, this sum has at most $N$ non-zero terms). Since $p(x,y)\geq 0$ for all $x,y\in \Gamma$ and (\ref{stochastic}) holds, one has, for all $p\in [1,+\infty]$ and all $f\in L^p(\Gamma)$,
\begin{equation}\label{contraction}
\left\Vert Pf\right\Vert_{L^p(\Gamma)}\leq \left\Vert f\right\Vert_{L^p(\Gamma)}.
\end{equation} 
We make use of the operator $P$ to define a Laplacian on $\Gamma$. Consider a function $f\in L^2(\Gamma)$. By (\ref{contraction}), $(I-P)f\in L^2(\Gamma)$ and
\begin{equation} \label{byparts}
\begin{array}{lll}
\displaystyle \langle (I-P)f,f\rangle_{L^2(\Gamma)} & = & \displaystyle \sum_{x,y} p(x,y)(f(x)-f(y))f(x)m(x)\\
& = & \displaystyle \frac 12 \sum_{x,y} p(x,y)\left\vert 
f(x)-f(y)\right\vert^{2} m(x),
\end{array}
\end{equation}
where we use (\ref{stochastic}) in the first equality and 
(\ref{reversibility}) in the second one. If we define now the operator ``length of the gradient''  by
\[
\nabla f(x)=\left(\frac 12\sum\limits_{y\in \Gamma} p(x,y)\left\vert 
f(y)-f(x)\right\vert^{2}\right)^{1/2}
\] 
for all function $f$ on $\Gamma$ and all $x\in \Gamma$ (this definition is taken from \cite{cougri}), (\ref{byparts}) shows that
\begin{equation} \label{bypartsbis}
\langle (I-P)f,f\rangle_{L^2(\Gamma)}=\left\Vert \nabla f\right\Vert_{L^2(\Gamma)}^2.
\end{equation}
Because of (\ref{reversibility}), the operator $P$ is self-adjoint on 
$L^2(\Gamma)$ and $I-P$, which, by (\ref{bypartsbis}) , can be considered as a discrete 
``Laplace'' operator, is non-negative and self-adjoint on $L^2(\Gamma)$. By means of spectral 
theory, one defines its square root $(I-P)^{1/2}$. The equality (\ref{bypartsbis}) exactly means that
\begin{equation} \label{sobl2}
\left\Vert (I-P)^{1/2}f\right\Vert_{L^2(\Gamma)}=\left\Vert \nabla f\right\Vert_{L^2(\Gamma)}.
\end{equation}
This equality has an interpretation in terms of Sobolev spaces defined through $\nabla$. Let $1\leq p\leq +\infty$. Say that a scalar-valued function $f$ on $\Gamma$ belongs to the (inhomogeneous) Sobolev space $W^{1,p}(\Gamma)$ (see also \cite{ostr}, \cite{gt}) if and only if
\[
\left\Vert f\right\Vert_{W^{1,p}(\Gamma)}:=\left\Vert f\right\Vert_{L^p(\Gamma)}+\left\Vert \nabla f\right\Vert_{L^p(\Gamma)} <+\infty.
\]
If $B$ is any ball in $\Gamma$ and $1\leq p\leq +\infty$, denote by $W^{1,p}_0(B)$ the subspace of $W^{1,p}(\Gamma)$ made of functions supported in $B$.\par
\noindent We will also consider the homogeneous versions of Sobolev spaces. For $1\leq p\leq +\infty$, define $\dot{E}^{1,p}(\Gamma)$ as the space of all scalar-valued functions $f$ on $\Gamma$ such that $\nabla f\in L^p(\Gamma)$, equipped with the semi-norm
\[
\left\Vert f\right\Vert_{\dot{E}^{1,p}(\Gamma)}:=\left\Vert \nabla f\right\Vert_{L^p(\Gamma)}.
\]
Then $\dot{W}^{1,p}(\Gamma)$ is the quotient space $\dot{E}^{1,p}(\Gamma)/\R$, equipped with the corresponding norm. 
%Note that, for $1\leq p<+\infty$, $\R\cap L^p(\Gamma)$ is equal to $0$ if $m(\Gamma)=\infty$ and to $\R$ if $m(\Gamma)<+\infty$, while $\R\cap L^{\infty}(\Gamma)=\R$. 
It is then routine to check that both inhomogeneous and homogeneous Sobolev spaces on $\Gamma$ are Banach spaces. \par
The equality (\ref{sobl2}) means that $\left\Vert (I-P)^{1/2}f\right\Vert_{L^2(\Gamma)}=\left\Vert f\right\Vert_{\dot{E}^{1,2}(\Gamma)}$. In other words, for $p=2$, the Sobolev spaces defined by $\nabla$ and by the Laplacian coincide. In the sequel, we address the analogous question for $p\neq 2$.
\subsection{Statement of the problem}
To that purpose, we consider separately two inequalities, the validity of which 
will be discussed in the sequel. Let $1<p<+\infty$. The first inequality we look at says that there 
exists $C_{p}>0$ such that, for all function $f$ on $\Gamma$ such that $(I-P)^{1/2}f\in L^p(\Gamma)$,
\begin{equation} \label{rp} \tag{$R_{p}$}
\left\Vert \nabla f\right\Vert_{p}\leq C_{p}\left\Vert
(I-P)^{1/2}f\right\Vert_{p}.
\end{equation}
This inequality means that the operator 
$\nabla (I-P)^{-1/2}$, which is nothing but the {\it Riesz transform} 
associated with $(I-P)$, is $L^p(\Gamma)$-bounded. Here and after, say that a (sub)linear 
operator $T$ is $L^p$-bounded, or is of strong type $(p,p)$, if there exists $C>0$ such that $\left\Vert Tf\right\Vert_p\leq C\left\Vert f\right\Vert_p$ for 
all $f\in L^p(\Gamma)$. Say that it is of weak type $(p,p)$ if there 
exists $C>0$ such that $m\left(\left\{x\in \Gamma,\ \left\vert 
Tf(x)\right\vert>\lambda\right\}\right)\leq \frac C{\lambda^p}\left\Vert f\right\Vert_{p}^p$ 
for all $f\in L^p(\Gamma)$ and all $\lambda>0$. Notice that he functions $f$ will be defined on $\Gamma$, whereas $Tf$ may be defined on $\Gamma$ or on $E$. \par
\noindent The second inequality under consideration says that there 
exists $C_{p}>0$ such that, for all function $f\in \dot{E}^{1,p}(\Gamma)$,
\begin{equation} \label{rrp} \tag{$RR_{p}$}
\left\Vert
(I-P)^{1/2}f\right\Vert_{p}\leq C_{p}\left\Vert \nabla f\right\Vert_{p}.
\end{equation}
(The notations $(R_p)$ and $(RR_p)$ are borrowed from \cite{ac}.) We 
have just seen, by (\ref{sobl2}), that $(R_{2})$ and $(RR_{2})$ always hold. A 
well-known fact (see \cite{r} for a proof in this context) is that, 
if (\ref{rp}) holds for some $1<p<+\infty$, then (\ref{rrp}) 
holds with $p^{\prime}$ such that $1/p+1/p^{\prime}=1$, while the converse is unclear in this discrete situation (it is false in the case of Riemannian manifolds, see \cite{ac}). As we will see, we have to consider four distinct issues: (\ref{rp}) for $p<2$, (\ref{rp}) for $p>2$, (\ref{rrp}) for $p<2$, (\ref{rrp}) for $p>2$.
\subsection{The $L^p$-boundedness of the Riesz transform}
\subsubsection{The case when $p<2$}
Let us first consider (\ref{rp}) when $p<2$. This problem was 
dealt with in \cite{r}, and we just recall the result proved therein, which involves 
some further assumptions on $\Gamma$. The first one is of geometric 
nature. Say that $(\Gamma,\mu)$ satisfies the doubling property if there exists $C>0$ such that, for all $x\in \Gamma$ and all $r>0$,  
\begin{equation} \label{D} \tag{$D$}
V(x,2r)\leq CV(x,r).
\end{equation}
Note that this assumption implies that there exist $C,D>0$ such that, for all $x\in \Gamma$, all $r>0$ and all $\theta>1$,
\begin{equation} \label{Dbis}
V(x,\theta r)\leq C\theta^DV(x,r).
\end{equation}
\begin{rem} \label{infinite}
Observe also that, since $(\Gamma,\mu)$ is infinite, it is also unbounded (since it is locally uniformly finite) so that, if (\ref{D}) holds, then $m(\Gamma)=+\infty$ (see \cite{martellthesis}). 
\end{rem}
The second assumption on $(\Gamma,\mu)$ is a uniform lower bound for 
$p(x,y)$ when $x\sim y$, {\it i.e.} when $p(x,y)>0$. For $\alpha>0$, say that $(\Gamma,\mu)$ satisfies the condition $\Delta(\alpha)$ if, for all $x,y\in \Gamma$,
\begin{equation} \label{deltaalpha} \tag{$\Delta(\alpha)$}
\left(x\sim y\Leftrightarrow \mu_{xy}\geq \alpha m(x)\right)\mbox{ and }x\sim x.
\end{equation}
The next two assumptions on $(\Gamma,\mu)$ are pointwise upper bounds 
for the iterates of $p$. Say that $(\Gamma,\mu)$ satisfies $(DUE)$ (a on-diagonal upper estimate for the iterates of $p$) if there exists $C>0$ such that, for all $x\in \Gamma$ and all $k\in \N^{\ast}$,
\begin{equation} \label{diagupper} \tag{$DUE$}
p_k(x,x)\leq \frac {Cm(x)}{V(x,\sqrt{k})}.
\end{equation}
Say that $(\Gamma,\mu)$ satisfies $(UE)$ (an upper estimate for the iterates of $p$) if there exist 
$C,c>0$ such that, for all $x,y\in \Gamma$ and all $k\in \N^{\ast}$,
\begin{equation} \label{upper} \tag{$UE$}
p_k(x,y)\leq \frac {Cm(x)}{V(x,\sqrt{k})}e^{-c\frac{d^{2}(x,y)}k}.
\end{equation}
Recall that, under assumption (\ref{D}), estimates (\ref{diagupper}) 
and (\ref{upper}) are equivalent (and the conjunction of (\ref{D}) and (\ref{diagupper}) is also equivalent to a Faber-Krahn inequality, \cite{cougri}, Theorem 1.1). The 
following result holds:
\begin{theo} \label{rieszpleq2} (\cite{r})
Under assumptions (\ref{D}), 
(\ref{deltaalpha}) and (\ref{diagupper}), (\ref{rp}) holds for all $1<p\leq 
2$. 
Moreover, the Riesz transform is of weak $(1,1)$ type, which means 
that there exists $C>0$ such that, for all $\lambda>0$ and all 
function $f\in L^{1}(\Gamma)$,
\[
m\left(\left\{x\in \Gamma;\ \nabla
(I-P)^{-1/2}f(x)>\lambda\right\}\right)\leq \frac C{\lambda} 
\left\Vert f\right\Vert_{1}.
\]
As a consequence, under the same assumptions, (\ref{rrp}) holds for 
all $2\leq p<+\infty$.
\end{theo}
Notice that, according to \cite{hebsal}, the assumptions of Theorem \ref{rieszpleq2} hold, for instance, when $\Gamma$ is the Cayley graph of a group with polynomial volume growth and $p(x,y)=h(y^{-1}x)$, where $h$ is a symmetric bounded probability density supported in a ball and bounded from below by a positive constant on an open generating neighborhood of $e$, the identity element of $G$, and actually Theorem \ref{rieszpleq2} had already been proved in \cite{hebsal}.
\subsubsection{The case when $p>2$}
When $p>2$, assumptions (\ref{D}), (\ref{upper}) and (\ref{deltaalpha}) are not 
sufficient to ensure the validity of (\ref{rp}), as the example of two 
copies of $\Z^{2}$ linked between with an edge shows (see \cite{r}, 
Section 4). More 
precisely, in this example, as explained in Section 4 of \cite{r}, the validity of (\ref{rp}) for $p>2$ would imply an $L^{2}$ Poincar\'e inequality on 
balls. \par
\noindent Say that $(\Gamma,\mu)$ satisfies a scaled $L^{2}$ Poincar\'e 
inequality on balls (this inequality will be denoted by $(P_{2})$ in 
the sequel) if there exists 
$C>0$ such that, for any $x\in \Gamma$, any $r>0$ and any function 
$f$ locally square integrable on $\Gamma$ such that $\nabla f$ is 
locally square integrable on $E$,
\begin{equation} \label{poincarel2} \tag{$P_{2}$}
\sum_{y\in B(x,r)} \left\vert f(y)-f_B\right\vert^{2}m(y)\leq 
Cr^{2}\sum_{y\in B(x,r)} \left\vert \nabla f(y)\right\vert^{2}m(y),
\end{equation}
where
\[
f_B=\frac 1{V(B)} \sum_{x\in B} f(x)m(x)
\]
is the mean value of $f$ on $B$. Under assumptions (\ref{D}), (\ref{poincarel2}) and (\ref{deltaalpha}), not only does 
(\ref{upper}) hold, but the iterates of $p$ also satisfy a pointwise Gaussian lower 
bound. Namely, there exist 
$c_{1},C_{1},c_{2},C_{2}>0$ such that, for all $n\geq 1$ and all 
$x,y\in \Gamma$ with $d(x,y)\leq n$,
\begin{equation} \label{lowerupper} \tag{$LUE$}
\frac {c_{1}m(x)}{V(x,\sqrt{n})}e^{-C_{1}\frac{d^{2}(x,y)}n}\leq 
p_{n}(x,y)\leq \frac {C_{2}m(x)}{V(x,\sqrt{n})}e^{-c_{2}\frac{d^{2}(x,y)}n}.
\end{equation}
Actually, (\ref{lowerupper}) is equivalent to the conjunction of 
(\ref{D}), (\ref{poincarel2}) and (\ref{deltaalpha}), and also to a 
discrete parabolic Harnack inequality, see \cite{delmo} (see also 
\cite{auschcou} for another approach of (\ref{lowerupper})). \par
Let $p>2$ and assume that (\ref{rp}) holds. Then, if $f\in 
L^p(\Gamma)$ and $n\geq 1$,
\begin{equation} \label{condness} \tag{$G_{p}$}
\left\Vert \nabla P^nf\right\Vert_{p}\leq \frac {C_p}{\sqrt{n}} \left\Vert 
f\right\Vert_{p}.
\end{equation}
Indeed, (\ref{rp}) implies that
\[
\left\Vert \nabla P^nf\right\Vert_{p}\leq C_p\left\Vert 
(I-P)^{1/2}P^nf\right\Vert_{p},
\]
and, due to the analyticity of 
$P$ on $L^p(\Gamma)$, one also has
\[
\left\Vert 
(I-P)^{1/2}P^nf\right\Vert_{p}\leq \frac {C'_p}{\sqrt{n}} \left\Vert 
f\right\Vert_{p}.
\]
More precisely, as was explained in \cite{r}, assumption $\Delta(\alpha)$ implies 
that $-1$ does not belong to the spectrum of $P$ on $L^{2}(\Gamma)$. 
As a consequence, $P$ is analytic on $L^{2}(\Gamma)$ (see 
\cite{cs}, Proposition 3), and since $P$ is submarkovian, $P$ is also 
analytic on $L^p(\Gamma)$ (see \cite{cs}, p. 426). Proposition 2 in 
\cite{cs} therefore yields the second inequality in (\ref{condness}). Thus, condition (\ref{condness}) is necessary for (\ref{rp}) to hold. Our first result is that, under assumptions (\ref{D}), 
(\ref{poincarel2}) and (\ref{deltaalpha}), for all $q>2$, condition $(G_q)$ is also sufficient 
for $(R_{p})$ to hold for all $2<p<q$:
\begin{theo} \label{rieszpgeq2}
Let $p_{0}\in (2,+\infty]$. Assume that $(\Gamma,\mu)$ satisfies (\ref{D}), 
(\ref{poincarel2}), (\ref{deltaalpha}) and $(G_{p_{0}})$. Then, for all 
$2\leq p<p_{0}$, 
(\ref{rp}) holds. As a consequence, if $p_{0}^{\prime}$ is such that 
$1/p_{0}+1/p_{0}^{\prime}=1$, (\ref{rrp}) holds for all 
$p_{0}^{\prime}<p\leq 2$.
\end{theo}
An immediate consequence of Theorem \ref{rieszpgeq2} and the previous discussion is the following result:
\begin{theo} \label{rieszpgeq2equiv}
Assume that $(\Gamma,\mu)$ satisfies (\ref{D}), 
(\ref{poincarel2}) and (\ref{deltaalpha}). Let $p_0\in (2,+\infty]$. Then, the following two assertions are equivalent:
\begin{itemize}
\item[$(i)$]
for all $p\in (2,p_0)$, $(G_p)$ holds,
\item[$(ii)$]
for all $p\in (2,p_0)$, $(R_p)$ holds.
\end{itemize}
\end{theo}
\begin{rem}
In the recent work \cite{dungeylp}, property $(G_p)$ is shown to be true for all $p\in (1,2]$ under the sole assumption that $\Gamma$ satisfies a local doubling property for the volume of balls.
\end{rem}
\begin{rem}
On Riemannian manifolds, the $L^2$ Poincar\'e inequality on balls is neither necessary, nor sufficient to ensure that the Riesz transform is $L^p$-bounded for all $p\in (2,\infty)$, see \cite{ac} and the references therein. We do not know if the corresponding assertion holds in the context of graphs.
\end{rem}
\subsubsection{Riesz transforms and harmonic functions}
We also obtain another characterization of the validity of $(R_p)$ for $p>2$ in terms of {\bf reverse H\"older inequalities for the gradient of harmonic functions}, in the spirit of \cite{shen} (in the Euclidean context for second order elliptic operators in divergence form) and \cite{ac} (on Riemannian manifolds). If $B$ is any ball in $\Gamma$ and $u$ a function on $B$, say that $u$ is harmonic on $B$ if, for all $x\in B\setminus \partial B$,
\begin{equation} \label{harmonic}
(I-P)u(x)=0.
\end{equation}
We will prove the following result:
\begin{theo} \label{rpharmonic}
Assume that (\ref{D}), (\ref{deltaalpha}) and (\ref{poincarel2}) hold. Then, there exists $p_0\in (2,+\infty]$ such that, for all $q\in (2,p_0)$, the following two conditions are equivalent:
\begin{itemize}
\item[$1.$]
$(R_p)$ holds for all $p\in (2,q)$,
%\item[$2.$]
%$(\Pi_p)$ holds for all $p\in (2,q)$,
\item[$2.$]
for all $p\in (2,q)$, there exists $C_p>0$ such that, for all ball $B\subset \Gamma$, all function $u$ harmonic in $32B$,
\begin{equation} \label{reverseholder} \tag{$RH_p$}
\left(\frac 1{V(B)} \sum_{x\in B} \left\vert \nabla u(x)\right\vert^p m(x)\right)^{\frac 1p} \leq C_p\left(\frac 1{V(16B)} \sum_{x\in 16B} \left\vert \nabla u(x)\right\vert^2 m(x)\right)^{\frac 12}.
\end{equation}
\end{itemize}
\end{theo}
Assertion $3.$ says that the gradient of any harmonic function in $32B$ satisfies a reverse H\"older inequality. Remember that such an inequality always holds for solutions of $\mbox{div}(A\nabla u)=0$ on any ball $B\subset\R^n$, if $u$ is assumed to be in $H^1(B)$ and $A$ is bounded and uniformly elliptic (see \cite{meyers}). In the present context, a similar self-improvement result can be shown:
\begin{pro} \label{selfimprove}
Assume that (\ref{D}), (\ref{deltaalpha}) and (\ref{poincarel2}) hold. Then there exists $p_0>2$ such that (\ref{reverseholder}) holds for any $p\in (2,p_0)$. As a consequence, (\ref{rp}) holds for any $p\in (2,p_0)$.
\end{pro}
As a corollary of Theorem \ref{rieszpleq2} and Proposition \ref{selfimprove}, we get:
\begin{cor} \label{reverserieszeps}
Assume that (\ref{D}), (\ref{deltaalpha}) and $(P_2)$ hold. Then, 
there exists $\varepsilon>0$ such that, for all $2-\varepsilon<p< 2+\varepsilon$, $\left\Vert \nabla f\right\Vert_{p}\sim \left\Vert 
(I-P)^{1/2}f\right\Vert_{p}$.
\end{cor}
\subsection{The reverse inequality}
Let us now focus on (\ref{rrp}). As already seen, (\ref{rrp}) holds 
for all $p>2$ under (\ref{D}), (\ref{deltaalpha}) and 
(\ref{diagupper}), and for all $p_{0}^{\prime}<p<2$ under 
(\ref{D}), (\ref{poincarel2}), (\ref{deltaalpha}) and ($G_{p_{0}}$) 
if $p_{0}>2$ and $1/p_{0}+1/p_{0}^{\prime}=1$. However, we can also 
give a sufficient condition for (\ref{rrp}) to hold for all $p\in (q_0,2)$ (for some $q_0<2$) which does not 
involve any assumption such that $(G_{p_{0}})$. For $1\leq p<+\infty$, 
say that $(\Gamma,\mu)$ satisfies a scaled $L^p$ Poincar\'e inequality on 
balls (this inequality will be denoted by $(P_{p})$ in 
the sequel) if there exists 
$C>0$ such that, for any $x\in \Gamma$, any $r>0$ and any function 
$f$ on $\Gamma$ such that $\left\vert f\right\vert^p$ and $\left\vert \nabla f\right\vert^p$ are locally integrable on $\Gamma$,
\begin{equation} \label{poincarelp} \tag{$P_{p}$}
\sum_{y\in B(x,r)} \left\vert f(y)-f_B\right\vert^{p}m(y)\leq 
Cr^{p}\sum_{y\in B(x,r)} \left\vert \nabla f(y)\right\vert^{p}m(y).
\end{equation}  
If $1\leq p<q<+\infty$, 
then $(P_{p})$ implies $(P_{q})$ (this is a very general statement on 
spaces of homogeneous type, {\it i.e.} on metric measured spaces where (\ref{D}) holds, see \cite{hk}). The converse implication 
does not hold but an $L^p$ Poincar\'e inequality still has a 
self-improvement in the following sense:
\begin{pro} \label{poincareself}
Let $(\Gamma,\mu)$ satisfy (\ref{D}). Then, for all $p\in 
(1,+\infty)$, if $(P_{p})$ holds, there exists $\varepsilon>0$ such 
that $(P_{p-\varepsilon})$ holds.
\end{pro}
This deep result actually holds in the general context of spaces of 
homogeneous type, {\it i.e.} when (\ref{D}) holds, see \cite{kz}. \par
Assuming that $(P_{q})$ holds for some $q<2$, we establish $(RR_{p})$ 
for $q<p<2$:
\begin{theo} \label{reverserieszpleq2}
Let $1\leq q<2$. Assume that (\ref{D}), (\ref{deltaalpha}) and $(P_q)$ hold. Then, for all $q<p<2$, $(RR_p)$ holds. Moreover, there exists $C>0$ such that, for all $\lambda>0$,
\begin{equation} \label{enoughestimate}
m\left(\left\{x\in \Gamma;\ \left\vert (I-P)^{1/2}f(x)\right\vert>\lambda\right\}\right)\leq \frac C{\lambda^q}\left\Vert \nabla f\right\Vert_{q}^q.
\end{equation}
\end{theo}
As a corollary of Theorem \ref{rieszpleq2}, Proposition \ref{poincareself} and Theorem \ref{reverserieszpleq2}, we get the following consequence:
\begin{cor} \label{reverserieszepsbis}
Assume that (\ref{D}), (\ref{deltaalpha}) and $(P_p)$ hold for some $p\in (1,2)$. Then, 
there exists $\varepsilon>0$ such that, for all $p-\varepsilon<q<+\infty$, $(RR_q)$ holds. In particular, $(RR_p)$ holds.
\end{cor}
%In this corollary, the notation $\left\Vert \nabla f\right\Vert_{p}\sim \left\Vert 
%(I-P)^{1/2}f\right\Vert_{p}$ means that there exists $C_p>0$ such that
%\[
%C_p^{-1}\left\Vert df\right\Vert_p\leq \left\Vert (I-P)^{1/2}f\right\Vert_p\leq C_p\left\Vert df\right\Vert_p
%\]
%for all function $f$.\par

%The following property of the gradient will be useful in the sequel. If $f,g$ are two functions on $\Gamma$ and $x\in \Gamma$, then
%\begin{equation} \label{derivation}
%\nabla (fg)(x)\leq \left\vert f(x)\right\vert \nabla g(x)+\sup_{y\sim x} \left\vert g(y)\right\vert \nabla f(x).
%\end{equation}
\subsection{An overview of the method}
Let us briefly describe the proofs of our results. Let us first 
consider Theorem \ref{rieszpgeq2}. The operator $T=\nabla(I-P)^{-1/2}$ 
can be written as
\[
T=\nabla \left(\sum_{k=0}^{+\infty} a_{k}P^k\right),
\]
where the $a_{k}$'s are defined by the expansion
\begin{equation} \label{expansion}
(1-x)^{-1/2}=\sum_{k=0}^{+\infty} a_{k}x^k
\end{equation}
for $-1<x<1$. The kernel of $T$ is therefore given by
\[
\nabla_x \left(\sum_{k=0}^{+\infty} a_{k}p_{k}(x,y)\right).
\]
It was proved in \cite{pota} that, under (\ref{D}) and (\ref{poincarel2}), this kernel satisfies the H\"ormander integral condition, which implies the $H^1(\Gamma)-L^1(\Gamma)$ boundedness of $T$ and therefore its $L^p(\Gamma)$-boundedness for all $1<p<2$, where $H^1(\Gamma)$ denotes the Hardy space on $\Gamma$ defined in the sense of Coifman and Weiss (\cite{coifman1}). However, the H\"ormander integral condition does not yield any information on the $L^p$-boundedness of $T$ for $p>2$. The proof of Theorem \ref{rieszpgeq2} actually relies on a theorem due to Auscher, Coulhon, Duong and Hofmann (\cite{acdh}), which, given some $p_0\in (2,+\infty]$, provides sufficient conditions for an $L^2$-bounded sublinear operator to be $L^p$-bounded for $2<p<p_0$. Let us recall this theorem here in the form to be used in the sequel for the sake of completeness (see \cite{acdh}, Theorem 2.1, \cite{auscher}, Theorem 2.2):
\begin{theo} \label{CZgeq2}
Let $p_0\in \left(2,+\infty\right]$. Assume that $\Gamma$ satisfies the doubling property (\ref{D}) and let $T$ be a sublinear operator acting on $L^2(\Gamma)$. For any ball $B$, let $A_B$ be a linear operator acting on $L^2(\Gamma)$, and assume that there exists $C>0$ such that, for all $f\in L^2(\Gamma)$, all $x\in \Gamma$ and all ball $B\ni x$,
\begin{equation} \label{cond1geq2}
\frac 1{V^{1/2}(B)} \left\Vert T(I-A_B)f\right\Vert_{L^2(B)}\leq C\left({\mathcal M}(\left\vert f\right\vert^2)\right)^{1/2}(x)
\end{equation}
and
\begin{equation} \label{cond2geq2}
\frac 1{V^{1/p_0}(B)} \left\Vert TA_Bf\right\Vert_{L^{p_0}(B)}\leq C\left({\mathcal M}(\left\vert Tf\right\vert^2)\right)^{1/2}(x).
\end{equation}
If $2<p<p_0$ and if, for all $f\in L^p(\Gamma)$, $Tf\in L^p(\Gamma)$, then there exists $C_p>0$ such that, for all $f\in L^2(\Gamma)\cap L^p(\Gamma)$,
\[
\left\Vert Tf\right\Vert_{L^p(\Gamma)}\leq C_p\left\Vert f\right\Vert_{L^p(\Gamma)}.
\]
\end{theo}
Notice that, to simplify the notations in our foregoing proofs, the formulation of Theorem \ref{CZgeq2} is slightly different from the one given in \cite{auscher} and in \cite{acdh}, since the family of operators $(A_r)_{r>0}$ used in these papers is replaced by a family $(A_B)$ indexed by the balls $B\subset \Gamma$, see Remark 5 after Theorem 2.2 in \cite{auscher}. Observe also that this theorem extends to vector-valued functions (this will be used in Section \ref{LPaley}). Finally, here and after, ${\mathcal M}$ denotes the Hardy-Littlewood maximal function: for any locally integrable function $f$ on $\Gamma$ and any $x\in \Gamma$,
\[
{\mathcal M}f(x)=\sup_{B\ni x} \frac 1{V(B)}\sum_{y\in B} \left\vert f(y)\right\vert m(y),
\]
where the supremum is taken over all balls $B$ containing $x$. Recall that, by the Hardy-Littlewood maximal theorem, since (\ref{D}) holds, ${\mathcal M}$ is of weak type $(1,1)$ and of strong type $(p,p)$ for all $1<p\leq +\infty$. \par
\noindent Following the proof of Theorem 2.1 in \cite{acdh}, we will obtain Theorem \ref{rieszpgeq2}  by applying Theorem \ref{CZgeq2} with $A_B=I-(I-P^{k^2})^n$ 
where $k$ is the radius of $B$ and $n$ is an integer only depending from the constant $D$ in (\ref{Dbis}).  \par
As far as Theorem \ref{reverserieszpleq2} is concerned, note first that $(RR_p)$ cannot be derived from $(R_{p^{\prime}})$ in this situation (where $1/p+1/p^{\prime}=1$), since we do not know whether $\left(R_{p^{\prime}}\right)$ holds or not under these assumptions. Following \cite{ac}, we first prove (\ref{enoughestimate}). The proof relies on a Calder\'on-Zygmund decomposition for Sobolev functions, which is the adaptation to our context of Proposition 1.1 in \cite{ac} (see also \cite{a} in the Euclidean case and \cite{am} for the extension to a weighted Lebesgue measure):
\begin{pro} \label{CZ}
Assume that (\ref{D}) and $(P_{q})$ hold for some $q\in [1,\infty)$ and let $p\in [q,+\infty)$. Let  $f \in \dot{E}^{1,p}(\Gamma)$ and $\alpha>0$. Then one can find a collection of balls $(B_{i})_{i\in I}$, functions $(b_{i})_{i\in I}\in \dot{E}^{1,q}(\Gamma)$ and a function $g\in\dot{E}^{1,\infty}$ such that the following properties hold:
\begin{equation}
f = g+\sum_{i\in I}b_{i}, \label{df}
\end{equation}
\begin{equation}
\|\nabla g\|_{\infty}\leq C\alpha, \label{eg}
\end{equation}
\begin{equation}
\supp b_{i}\subset B_{i}, \,\sum_{x\in B_{i}}| \nabla b_{i}|^{q}(x)m(x)\leq C\alpha^{q}V(B_{i}),\label{eb}
\end{equation}
\begin{equation}
\sum_{i\in I}V(B_{i})\leq C\alpha^{-p}\sum_{x\in \Gamma} | \nabla f|^{p}(x) m(x),
\label{eB}
\end{equation}
\begin{equation}
\sum_{i\in I}\chi_{B_{i}}\leq N, \label{rb}
\end{equation}
where $C$ and $N$  only depend on $q$, $p$ and on the constants in $(D)$ and $(P_{q})$.
\end{pro}
As in \cite{ac}, we rely on this Calder\'on-Zygmund decomposition to 
establish (\ref{enoughestimate}). The argument also uses the $L^p(\Gamma)$-boundedness, for all $2<p<+\infty$, of a discrete version of the Littlewood-Paley-Stein $g$-function (see \cite{topics}), 
which does not seem to have been stated before in this context and is interesting in itself. For all function $f$ on $\Gamma$ and all $x\in \Gamma$, define
\[
g(f)(x)=\left(\sum_{l\geq 1} l\left\vert (I-P)P^lf(x)\right\vert^2\right)^{1/2}.
\]
Observe that this is indeed a discrete analogue of the $g$-function 
introduced by Stein in \cite{topics}, since $(I-P)P^l=P^l-P^{l+1}$ 
can be seen as a discrete time derivative of $P^l$ and $P$ is a Markovian operator. \par
\noindent It is easy to check that the sublinear operator $g$ is bounded 
in $L^2(\Gamma)$. Indeed, as  already said, the assumption (\ref{deltaalpha}) implies that the spectrum of $P$ is contained in $[a,1]$ for some $a>-1$. As 
a consequence, $P$ can be written as
\[
P=\int_a^1 \lambda dE(\lambda),
\]
so that, for all integer $l\geq 1$,
\[
(I-P)P^l=\int_a^1 (1-\lambda)\lambda^ldE(\lambda)
\]  
and, for all $f\in L^2(\Gamma)$,
\[
\left\Vert (I-P)P^lf\right\Vert_2^2=\int_a^1 (1-\lambda)^2\lambda^{2l} dE_{f,f}(\lambda).
\]
It follows that, for all $f\in L^2(\Gamma)$,
\[
\begin{array}{lll}
\left\Vert g(f)\right\Vert_2^2 & = & \displaystyle \sum_{l\geq 1} l \left\Vert (I-P)P^lf\right\Vert_2^2\\
& = & \displaystyle \int_a^1 (1-\lambda)^2\sum_{l\geq 1} l\lambda^{2l} dE_{f,f}(\lambda)\\
& = &\displaystyle \int_a^1 \left(\frac{\lambda}{1+\lambda}\right)^2dE_{f,f}(\lambda)\\
& \leq & \displaystyle \left\Vert f\right\Vert_2^2.
\end{array}
\]
It turns out that, as in the Littlewood-Paley-Stein semigroup theory, $g$ is also $L^p$-bounded for $1<p<+\infty$:
\begin{theo} \label{LP}
Assume that (\ref{D}), (\ref{diagupper}) and (\ref{deltaalpha}) 
hold. Let $1<p<+\infty$. There exists $C_p>0$ such that, for all $f\in L^p(\Gamma)$,
\[
\left\Vert g(f)\right\Vert_p\leq C_p \left\Vert f\right\Vert_p.
\]
\end{theo}
Actually, this inequality will only be used for $p>2$ in the sequel, but 
the result, which is interesting in itself, does hold and will be proved for all $1<p<+\infty$. \par
\noindent Before going further, let us mention that, in \cite{dungeylp}, N. Dungey establishes, under a local doubling property for the volume of balls, the $L^p$-boundedness for all $p\in (1,2]$ of another version of the Littlewood-Paley-Stein functional, involving the gradient instead of the ``time derivative'' and the (continuous time) semigroup generated by $I-P$.  Although we do not use Dungey's result here, it may prove useful to study the boundedness of Riesz transforms on graphs.\par
\noindent The proof of Theorem \ref{LP} for $p>2$ relies on the vector-valued version of Theorem \ref{CZgeq2}, while, for $p<2$, we use the vector-valued version of the following result (see \cite{auscher}, Theorem 2.1 and also \cite{bk} for an earlier version):
\begin{theo} \label{CZleq2}
Let $p_0\in [1,2)$. Assume that $\Gamma$ satisfies the doubling property (\ref{D}) and let $T$ be a sublinear operator of strong type $(2,2)$. For any ball $B$, 
let $A_B$ be a linear operator acting on $L^2(\Gamma)$. Assume that, 
for all $j\geq 1$, there exists $g(j)>0$ such that, for all ball $B\subset 
\Gamma$ and all function $f$ supported in $B$,
\begin{equation} \label{cond1leq2}
\frac 1{V^{1/2}(2^{j+1}B)} \left\Vert 
T(I-A_B)f\right\Vert_{L^2(C_j(B))}\leq g(j)\frac 1{V^{1/p_0}(B)} \left\Vert 
f\right\Vert_{L^{p_0}}
\end{equation}
for all $j\geq 2$ and
\begin{equation} \label{cond2leq2}
\frac 1{V^{1/2}(2^{j+1}B)}  \left\Vert A_Bf\right\Vert_{L^2(C_j(B))}\leq 
g(j)\frac 1{V^{1/p_0}(B)} \left\Vert f\right\Vert_{L^{p_{0}}}
\end{equation}
for all $j\geq 1$. If $\displaystyle \sum_{j\geq 1} g(j)2^{Dj}<+\infty$ where $D$ is given by (\ref{Dbis}), then $T$ is of weak type $(p_0,p_0)$, and is therefore of strong type $(p,p)$ for all $p_0<p<2$.
\end{theo}
Going back to Theorem \ref{reverserieszpleq2}, once 
(\ref{enoughestimate}) is established, we conclude by applying real interpolation theorems for Sobolev spaces, which are also new in this context. More precisely, we prove:
\begin{theo} \label{interpolation}
Let $q\in [1,+\infty)$ and assume that (\ref{D}), $(P_q)$ and (\ref{deltaalpha}) hold. Then, for all $q<p<+\infty$, $\dot{W}^{1,p}(\Gamma)=\left(\dot{W}^{1,q}(\Gamma),\dot{W}^{1,\infty}(\Gamma)\right)_{1-\frac qp,p}$.
%is an interpolation space between 
%$\dot{W}^{1,q}(\Gamma)$ and $\dot{W}^{1,\infty}(\Gamma)$.
\end{theo}
As an immediate corollary, we obtain:
\begin{cor}[The reiteration theorem]\label{RH} Assume that $\Gamma$ satisfies (\ref{D}), $(P_{q})$ for some $1\leq q<+\infty$ and (\ref{deltaalpha}). Define $q_{0}=\inf\left\lbrace q \in [1,\infty):(P_{q})\textrm{ holds}\right\rbrace$. For $q_{0}<p_{1}<p<p_{2}\leq +\infty$, if $\displaystyle \frac 1p=\frac{1-\theta}{p_1}+\frac{\theta}{p_2}$, then $\dot{W}^{1,p}(\Gamma)=\left(\dot{W}^{1,p_1}(\Gamma),\dot{W}^{1,p_2}(\Gamma)\right)_{\theta,p}$.
\end{cor}
Corollary \ref{RH}, in conjunction with (\ref{enoughestimate}), conclude the proof of Theorem \ref{reverserieszpleq2}. Notice that, since we know that Sobolev 
spaces interpolate by the real method, we do not need any argument as the one in Section 1.3 of \cite{ac}.\par
For the proof of Theorem \ref{rpharmonic}, we introduce a {\it discrete differential} and go through a property analogous to $(\Pi_p)$ in \cite{ac}, see Section \ref{rieszharmonic} for detailed definitions. Proposition \ref{selfimprove} follows essentially from Gehring's self-improvement of reverse H\"older inequalities (\cite{gehring}).\par
%Let us finally mention that, in \cite{ac}, in the context of Riemannian manifolds, the authors also obtain a characterization of the validity of $(R_p)$ for all $p\in (2,p_0)$ (with some $p_0>2$) in terms of a reverse H\"older inequality for harmonic functions on balls. We intend to discuss the corresponding issue on graphs in a forthcoming paper.\par
\noindent The plan of the paper is as follows. After recalling some well-known estimates for the iterates of $p$ and deriving some consequences (Section \ref{bounds}), we first prove Theorem \ref{LP}, 
which is of independent interest, in Section \ref{LPaley}. In Section \ref{riesz}, we prove Theorem \ref{rieszpgeq2} using Theorem \ref{CZgeq2}. Section \ref{decompo} is devoted to the proof of Proposition \ref{CZ}. Theorem \ref{interpolation} is established in Section \ref{interpo} by methods similar to \cite{badr1} and, in Section \ref{reverseriesz}, we prove Theorem \ref{reverserieszpleq2}. Finally, Section \ref{rieszharmonic} contains the proof of Theorem \ref{rpharmonic} and of Proposition \ref{selfimprove}.

\SE{Kernel bounds} \label{bounds}

In this section, we gather some estimates for the iterates of $p$ and some straightforward consequences of frequent use in the sequel. We always assume that (\ref{D}), (\ref{poincarel2}) and (\ref{deltaalpha}) hold. First, as already said, (\ref{lowerupper}) holds. Moreover, we also have the following pointwise estimate for the discrete ``time derivative'' of $p_l$: there exist $C,c>0$ such that, for all $x,y\in \Gamma$ and all $l\in \N^{\ast}$,
\begin{equation} \label{timederivative}
\left\vert p_l(x,y)-p_{l+1}(x,y)\right\vert \leq \frac{Cm(y)}{lV(x,\sqrt{l})} e^{-c\frac{d^2(x,y)}l}.
\end{equation}
This ``time regularity'' estimate, which is a consequence of the 
$L^{2}$ analyticity of $P$, was first proved by Christ 
(\cite{christ}) by a quite difficult argument. Simpler proofs have 
been given by Blunck (\cite{blunck}) and, more recently, by Dungey 
(\cite{dungey}). \par
Thus, if $B$ is a ball in $\Gamma$ with radius $k$, $f$ any function supported in $B$ and $i\geq 2$, one has, for all $x\in C_i(B)$ and all $l\geq 1$,
\begin{equation} \label{offdiagestim}
\left\vert P^lf(x)\right\vert + l\left\vert (I-P)P^lf(x)\right\vert\leq \frac C{V(B)}e^{-c\frac{4^ik^2}l}\left\Vert f\right\Vert_{L^1}.
\end{equation}
This ``off-diagonal'' estimate follows from (\ref{upper}) and (\ref{timederivative}) and the fact that, for all $y\in B$, by (\ref{D}),
\[
V(y,k)\sim V(B)\mbox{ and }\frac{V(y,k)}{V(y,\sqrt{l})}\leq C\sup\left(1,\left(\frac{k}{\sqrt{l}}\right)^D\right).
\]
Similarly, if $B$ is a ball in $\Gamma$ with radius $k$, $i\geq 2$ and $f$ any function supported in $C_i(B)$,  one has, for all $x\in B$ and all $l\geq 1$,
\begin{equation} \label{offdiagestimbis}
\left\vert P^lf(x)\right\vert + l\left\vert (I-P)P^lf(x)\right\vert\leq \frac C{V(2^{i}B)}e^{-c\frac{4^ik^2}l}\left\Vert f\right\Vert_{L^1}.
\end{equation}
Finally, for all ball $B$ with radius $k$, all $i\geq 2$, all function $f$ supported in $C_i(B)$ and all $l\geq 1$,
\begin{equation} \label{gradl2off}
\left\Vert  \nabla P^lf\right\Vert_{L^2(B)}\leq \frac C{\sqrt{l}} e^{-c\frac{4^{i}k^2}{l}}\left\Vert f\right\Vert_{L^2(C_i(B))}.
\end{equation}
See Lemma 2 in \cite{r}. If one furthermore assumes that $(G_{p_0})$ holds for some $p_0>2$, then, by interpolation between (\ref{gradl2off}) and $(G_{p_0})$, one obtains, for all $p\in (2,p_0)$, all $f$ supported in $C_i(B)$ and all $l\geq 1$,
\begin{equation} \label{gradlpoff}
\left\Vert  \nabla P^lf\right\Vert_{L^p(B)}\leq \frac {C_p}{\sqrt{l}} e^{-c\frac{4^{i}k^2}{l}}\left\Vert f\right\Vert_{L^p(C_i(B))}.
\end{equation}
Inequalities (\ref{gradl2off}) and (\ref{gradlpoff}) may be regarded as ``Gaffney'' type inequalities, in the spirit of \cite{gaffney}.

\SE{Littlewood-Paley inequalities} \label{LPaley}
In this section, we establish Theorem \ref{LP}.\par
\noindent {\bf The case $1<p<2$: } We apply the vector-valued version of 
Theorem \ref{CZleq2} with $T=g$ and $p_0=1$ and, for all ball $B$ with 
radius $k$, $A_{B}$ defined by
\[
A_B=I-(I-P^{k^2})^n,
\]
where $n$ is a positive integer, to be chosen in the proof. More 
precisely, we consider, for $f\in L^{2}(\Gamma)$ and $x\in \Gamma$,
\[
Tf(x)=\left(\sqrt{l}(I-P)P^lf(x)\right)_{l\geq 1},
\]
so that $T$ maps $L^{2}(\Gamma)$ into $L^{2}(\Gamma,l^{2})$.\par
\noindent Let $B$ be a ball and $f$ supported in $B$. Let us first check (\ref{cond1leq2}). Using the expansion
\[
(I-P^{k^2})^n=\sum_{p=0}^n C_n^p (-1)^pP^{pk^2},
\]
we obtain
\[
T(I-A_{B})f=\left(\alpha_l(I-P)P^lf\right)_{l\geq 1}
\]
where
\[
\alpha_l:=\sum_{0\leq p\leq n;\ l\geq pk^2} C_n^p (-1)^p\sqrt{l-pk^2}.
\]
Since it follows from (\ref{offdiagestim}) that
\begin{equation} \label{gf}
\frac 1{V(2^{j+1}B)} \left\Vert (I-P)P^lf\right\Vert_{L^2(C_j(B))}^2\leq \frac C{l^2V^2(B)}e^{-c\frac{4^jk^2}l}\left\Vert f\right\Vert_{L^1}^2,
\end{equation}
we will be able to go on thanks to the following estimate:
\begin{lem} \label{estimalphal}
There exists $C>0$ only depending on $n$ such that, for all $j\geq 2$,
\[
\sum_{l\geq 1} \frac{\left\vert \alpha_l\right\vert^2}{l^2}e^{-c\frac{4^jk^2}l}\leq C4^{-2nj}.
\]
\end{lem}
{\bf Proof of Lemma \ref{estimalphal}: } If $mk^2\leq l<(m+1)k^2$ for some integer $0\leq m\leq n$, one obviously has
\begin{equation} \label{alphal1}
\left\vert \alpha_l\right\vert\leq Ck\sqrt{m+1},
\end{equation}
where $C>0$ only depends on $n$, while, if $l>(n+1)k^2$, one has
\begin{equation} \label{alphal2}
\left\vert \alpha_l\right\vert\leq Cl^{-\frac{2n-1}2}k^{2n}.
\end{equation}
This estimate follows from the following inequality, valid for any $C^n$ function $\varphi$ on $(0,+\infty)$:
\begin{equation} \label{fellerineq}
\left\vert \sum_{p=0}^n C^p_n (-1)^p \varphi(t-pk^2)\right\vert\leq C\sup_{u\geq \frac t{n+1}} \left\vert \varphi^{(n)}(u)\right\vert k^{2n},
\end{equation}
where $C>0$ only depends on $n$ (see \cite{feller}, problem 16, p. 65). It follows from (\ref{alphal1}) that, for all $0\leq m\leq n$,
\[
\begin{array}{lll} 
\displaystyle \sum_{mk^2< l\leq (m+1)k^2}  \frac{\left\vert \alpha_l\right\vert^2}{l^2} e^{-c\frac{4^jk^2}l} &\leq & \displaystyle C\sum_{mk^2<l\leq (m+1)k^2} \frac{(m+1)k^2}{l^2} e^{-c\frac{4^jk^2}l} \\
& \leq & \displaystyle C\int_{mk^2}^{(m+1)k^2} \frac{(m+1)k^2}{t^2} e^{-c\frac{4^jk^2}t} dt\\
& \leq & \displaystyle  Ce^{-c4^j}
\end{array}
\]
where $C,c>0$ only depend on $n$. Similarly, thanks to (\ref{alphal2}),
\[
\begin{array}{lll}
\displaystyle \sum_{l>(n+1)k^2}  \frac{\left\vert \alpha_l\right\vert^2}{l^2} e^{-c\frac{4^jk^2}l} &\leq & \displaystyle C\sum_{l>(n+1)k^2} \frac{l^{-(2n-1)}k^{4n}}{l^2} e^{-c\frac{4^jk^2}l}\\
& \leq & \displaystyle C\int_{(n+1)k^2}^{+\infty} \frac{k^{4n}}{t^{2n+1}}e^{-c\frac{4^jk^2}t}dt\\
&\leq &\displaystyle C4^{-2nj} \int_0^{+\infty} \frac{e^{-c/w}}{w^{2n+1}}dw,\\
& = & \displaystyle C4^{-2nj},
\end{array}
\]
which concludes the proof of Lemma \ref{estimalphal}. \hfill\fin\par
\noindent Finally, one obtains
\[
\begin{array}{lll}
\displaystyle \frac 1{V^{1/2}(2^{j+1}B)} \left\Vert T(I-A_B)f\right\Vert_{L^2((2^{j+1}B\setminus 2^jB),l^2)} & \leq & \displaystyle C\frac{4^{-nj}}{V(B)}\left\Vert f\right\Vert_{L^1},\end{array}
\]
which means that  (\ref{cond1leq2}) holds with $g(j)=4^{-nj}$, and one just has to choose $\displaystyle n>\frac D2$ in order to have $\displaystyle \sum_j g(j)2^{Dj}<+\infty$.
\par
\noindent Let us now check (\ref{cond2leq2}). Since
\[
A_B=\sum_{p=1}^n C_n^p (-1)^pP^{pk^2},
\]
it is enough to prove that, for all $j\geq 1$ and all $1\leq p\leq n$,
\begin{equation} \label{enoughcond2leq2}
\frac 1{V^{1/2}(2^{j+1}B)}\left\Vert P^{pk^2}f\right\Vert_{L^2(C_j(B))}\leq g(j)\frac 1{V(B)}\left\Vert f\right\Vert_{L^1(B)}.
\end{equation}
For all $x\in C_j(B)$, (\ref{offdiagestim}) yields 
\[
\left\vert P^{pk^2}f(x)\right\vert \leq C\frac{e^{-c'\frac{4^j}p}}{V(B)}\left\Vert f\right\Vert_{L^1(B)}
\]
if $j\geq 2$, and
\[
\left\vert P^{pk^2}f(x)\right\vert \leq \frac C{V(B)} \left\Vert f\right\Vert_{L^1(B)}
\]
for $j=1$, just by (\ref{upper}). As a consequence,
\[
\left\Vert P^{pk^2}f\right\Vert_{L^2(C_j(B))}\leq C\frac{e^{-c'\frac{4^j}p}}{V(B)}V^{1/2}(2^{j+1}B)\left\Vert f\right\Vert_{L^1(B)},
\]
so that (\ref{enoughcond2leq2}) holds. This ends the proof of Theorem \ref{LP} when $1<p<2$.\par
{\bf The case $2<p<+\infty$: } This time, we apply the vector-valued version of Theorem \ref{CZgeq2} with the same choices of $T$ and $A_B$. Let us first check (\ref{cond1geq2}), which reads in this situation as
\[
\frac 1{V^{1/2}(B)}\left\Vert T(I-A_B)f\right\Vert_{L^2(B,l^2)}\leq C\left({\mathcal M}\left(\left\vert f\right\vert^2\right)\right)^{1/2}(y)
\]
for all $f\in L^2(\Gamma)$, all ball $B\subset \Gamma$ and all $y\in B$. Fix such an $f$, such a ball $B$ and $y\in B$. Write
\[
f=\sum_{j\geq 1} f\chi_{C_j(B)}:=\sum_{j\geq 1}f_j.
\]
The $L^2$-boundedness of $g$ and $A_B$ and the doubling property (\ref{D}) yield
\[
\frac 1{V^{1/2}(B)}\left\Vert T(I-A_B)f_1\right\Vert_{L^2(B,l^2)}\leq \frac C{V^{1/2}(B)}\left\Vert f\right\Vert_{L^2(4B)}\leq C\left({\mathcal M}\left(\left\vert f\right\vert^2\right)\right)^{1/2}(y).
\]
Let $j\geq 2$. Using the same notations as for the case $1<p<2$, one has
\[
\left\Vert T(I-A_B)f_j\right\Vert_{L^2(B,l^2)}^2=\sum_{l\geq 1} \left\vert \alpha_l\right\vert^2 \sum_{x\in B} \left\vert (I-P)P^lf_j(x)\right\vert^2m(x).
\]
For all $x\in B$, it follows from (\ref{offdiagestimbis}) and the Cauchy-Schwarz inequality that
\[
\begin{array}{lll}
\displaystyle \left\vert (I-P)P^lf_j(x)\right\vert & \leq & \displaystyle  \frac Cl e^{-c'\frac{4^jk^2}l} \frac 1{V(2^jB)} \sum_{z\in 2^{j+1}B} \left\vert f_j(z)\right\vert m(z)\\
& \leq & \displaystyle  \frac Cl e^{-c'\frac{4^jk^2}l} \frac 1{V^{1/2}(2^{j+1}B)} \left(\sum_{z\in 2^{j+1}B} \left\vert f_j(z)\right\vert^2 m(z)\right)^{1/2}\\
& \leq & \displaystyle \frac Cl e^{-c'\frac{4^jk^2}l} \ \left({\mathcal M} \left(\left\vert f\right\vert^2\right)(y)\right)^{1/2}.
\end{array}
\]
As a consequence, by Lemma \ref{estimalphal},
\[
\begin{array}{lll}
\displaystyle \left\Vert T(I-A_B)f_j\right\Vert_{L^2(B,l^2)}^2 &\leq & \displaystyle C \left(\sum_{l\geq 1} \frac{\left\vert \alpha_l\right\vert^2}{l^2}e^{-c\frac{4^jk^2}l}\right) {\mathcal M} \left(\left\vert f\right\vert^2\right)(y)V(B)\\
& \leq & \displaystyle CV(B)4^{-2nj}{\mathcal M} \left(\left\vert f\right\vert^2\right)(y),
\end{array}
\]
which yields (\ref{cond1geq2}) by summing up on $j\geq 1$.\par
\noindent To prove (\ref{cond2geq2}), it suffices to establish that, for all $1\leq j\leq n$, all ball $B\subset \Gamma$ and all $y\in B$,
\[
\left\Vert TP^{jk^2}f\right\Vert_{L^{\infty}(B,l^2)}\leq C\left({\mathcal M}\left\Vert Tf\right\Vert_{l^2}^2(y)\right)^{1/2}.
\]
Let $x\in B$. By Cauchy-Schwarz and the fact that 
\[
\sum_{y\in \Gamma} p_{jk^2}(x,y)=1
\]
for all $x\in \Gamma$, one has, for any function $h\in L^2(\Gamma)$,
\[
\left\vert P^{jk^2}h(x)\right\vert\leq \left(P^{jk^2}\left\vert h\right\vert^2(x)\right)^{1/2}.
\]
It follows that, for all $l\geq 1$,
\[
\left\vert P^{jk^2}(\sqrt{l}(I-P)P^lf)(x)\right\vert^2\leq P^{jk^2}\left(l\left\vert (I-P)P^lf\right\vert^2\right)(x),
\]
so that
\[
\begin{array}{lll}
\displaystyle \sum_{l\geq 1} \left\vert P^{jk^2}(\sqrt{l}(I-P)P^lf)(x)\right\vert^2 & \leq & \displaystyle P^{jk^2} \left(\sum_{l\geq 1} l\left\vert (I-P)P^lf\right\vert^2\right)(x)\\
& = & \displaystyle  P^{jk^2} \left(\left\Vert Tf\right\Vert_{l^2}^2\right)(x)\\
& \leq & \displaystyle C{\mathcal M}\left(\left\Vert Tf\right\Vert_{l^2}^2\right)(y),
\end{array}
\]
which is the desired estimate (note that the last inequality follows easily from (\ref{upper})). Thus, (\ref{cond2geq2}) holds and the proof of Theorem \ref{LP} is therefore complete. \hfill\fin

\SE{Riesz transforms for $p>2$} \label{riesz}
% Let $p\in (1,+\infty)$ and assume that $R$ is $L^p$-bounded. Then, for all $n\geq 1$ and all $f\in L^p(\Gamma)$,
% \[
% \left\Vert \left\vert \nabla P^nf\right\vert\right\Vert_{L^p(\Gamma)}\leq \left\Vert (I-P)^{1/2}P^nf\right\Vert_{L^p(\Gamma)}\leq \frac C{\sqrt{n}}\left\Vert f\right\Vert_{L^p(\Gamma)}.
% \]
In the present section, we establish Theorem \ref{rieszpgeq2}, applying Theorem \ref{CZgeq2} with the same choice of $A_B$ as in Section \ref{LPaley}. One has $\left\Vert A_B\right\Vert_{2,2}=1$. In view of Theorem \ref{CZgeq2}, it suffices to show that
\begin{equation} \label{cond1bis}
\frac 1{V^{1/2}(B)} \left\Vert T(I-P^{k^2})^nf\right\Vert_{L^2(B)}\leq C\left({\mathcal M}(\left\vert f\right\vert^2)\right)^{1/2}(x)
\end{equation}
and
\begin{equation} \label{cond2bis}
\frac 1{V^{1/p_0}(B)} \left\Vert T\left(I-(I-P^{k^2})^n\right)f\right\Vert_{L^{p_0}(B)}\leq C\left({\mathcal M}(\left\vert Tf\right\vert^2)\right)^{1/2}(x)
\end{equation}
for all $f\in L^2(\Gamma)$, all $x\in \Gamma$ and all ball $B\subset \Gamma$ containing $x$. Fix such data $f,x$ and $B$.\par
{\bf Proof of (\ref{cond1bis}): }Set $f_i=f\chi_{C_i(B)}$ for all $i\geq 1$. The $L^2$-boundedness of $T(I-P^{k^2})^n$ yields
\begin{equation} \label{cond1bis1}
\frac 1{V^{1/2}(B)} \left\Vert T(I-P^{k^2})^nf_1\right\Vert_{L^2(B)}\leq \frac C{V^{1/2}(B)} \left\Vert f_1\right\Vert_{L^2(\Gamma)}\leq C\left({\mathcal M}(\left\vert f\right\vert^2)\right)^{1/2}(x).
\end{equation}
Fix now $i\geq 2$. In order to estimate the left-hand side of (\ref{cond1bis}) with $f$ replaced by $f_i$, we use the expansion
\[
(I-P)^{-1/2}=\sum_{l=0}^{+\infty} a_lP^l,
\]
where the $a_l$'s are defined by (\ref{expansion}) (observe that, for all $l\geq 0$, $a_l>0$). Therefore, one has
\[
\begin{array}{lll}
\displaystyle (I-P)^{-1/2}(I-P^{k^2})^nf_i&= &\displaystyle \sum_{l=0}^{+\infty} a_l P^l(I-P^{k^2})^nf_i \\
& =& \displaystyle \sum_{l=0}^{+\infty} a_l\sum_{j=0}^n C^j_n (-1)^j P^{l+jk^2}f_i\\
& = & \displaystyle \sum_{l=0}^{+\infty} d_lP^lf_i,
\end{array}
\]
where
\[
d_l=\sum_{0\leq j\leq n,\ jk^2\leq l} (-1)^j C^j_n a_{l-jk^2}.
\]
It follows that
\[
\left\vert T(I-P^{k^2})^nf_i(x)\right\vert\leq \sum_{l=1}^{+\infty} \left\vert d_l\right\vert \nabla P^lf_i(x)
\]
for all $x\in B$. Indeed, if $x\in B$ and $l=0$, $\nabla P^lf_i(x)=\nabla f_i(x)=0$ because $f_i$ is supported in $C_i(B)$. Thus, one has
\[
\left\Vert T(I-P^{k^2})^nf_i\right\Vert_{L^2(B)}\leq \sum_{l=1}^{+\infty} \left\vert d_l\right\vert \left\Vert \left\vert \nabla P^lf_i\right\vert\right\Vert_{L^2(B)}.
\]
According to (\ref{gradl2off}), one has
\begin{equation} \label{cond1bis2}
\left\Vert T(I-P^{k^2})^nf_i\right\Vert_{L^2(B)}\leq C \sum_{l=1}^{+\infty} \left\vert d_l\right\vert  \frac{e^{-c\frac{4^{i}k^2}{l}}}{\sqrt{l}} \left\Vert f\right\Vert_{L^2(2^{i+1}B\setminus 2^{i}B)}.
\end{equation}
We claim that the following estimates hold for the $d_l$'s:
\begin{lem} \label{estim}
There exists $C>0$ only depending on $n$ with the following properties: for all integer $l\geq 1$,
\begin{itemize}
\item[$(i)$]
if there exists an integer $0\leq m\leq n$ such that $mk^2<l<(m+1)k^2$, $\left\vert d_l\right\vert\leq \frac C{\sqrt{l-mk^2}}$,
\item[$(ii)$]
if there exists an integer $0\leq m\leq n$ such that $l=(m+1)k^2$, $\left\vert d_l\right\vert \leq C$,
\item[$(iii)$]
if $l>(n+1)k^2$, $\left\vert d_l\right\vert\leq Ck^{2n}l^{-n-\frac 12}$.
\end{itemize}
\end{lem}
We postpone the proof of this lemma to the Appendix and end the proof of (\ref{cond1bis}). According to (\ref{cond1bis2}), one has
\begin{equation} \label{s1s2s3}
\begin{array}{lll}
\displaystyle \left\Vert T(I-P^{k^2})^nf_i\right\Vert_{L^2(B)} &\leq & \displaystyle C\sum_{m=0}^n \sum_{mk^2<l<(m+1)k^2} \left\vert d_l\right\vert  \frac{e^{-c\frac{4^{i}k^2}{l}}}{\sqrt{l}} \left\Vert f\right\Vert_{L^2(2^{i+1}B\setminus 2^{i}B)}\\
& + & \displaystyle C\sum_{m=0}^n \left\vert d_{(m+1)k^2}\right\vert \frac{e^{-c\frac{4^{i}}{m+1}}}{k\sqrt{m+1}} \left\Vert f\right\Vert_{L^2(2^{i+1}B\setminus 2^{i}B)}\\
& + & \displaystyle C\sum_{l>(n+1)k^2} \left\vert d_l\right\vert  \frac{e^{-c\frac{4^{i}k^2}{l}}}{\sqrt{l}} \left\Vert f\right\Vert_{L^2(2^{i+1}B\setminus 2^{i}B)}\\
& := & S_1+S_2+S_3.
\end{array}
\end{equation}
For $S_1$, Lemma \ref{estim} yields
\[
\left\vert S_1\right\vert  \leq C\sum_{m=0}^n \sum_{mk^2<l<(m+1)k^2} \frac {e^{-c\frac{4^{i}k^2}{l}}}{\sqrt{l}\sqrt{l-mk^2}}  \left\Vert f\right\Vert_{L^2(2^{i+1}B\setminus 2^{i}B)}.
\]
But, for each $1\leq m\leq n$, 
\[
\begin{array}{lll}
\displaystyle \sum_{mk^2<l<(m+1)k^2} \frac {e^{-c\frac{4^{i}k^2}{l}}}{\sqrt{l}\sqrt{l-mk^2}} & \leq &\displaystyle C\int_{mk^2}^{(m+1)k^2} \frac{e^{-c\frac{4^ik^2}t}}{\sqrt{t-mk^2}\sqrt{t}}dt\\
& \leq & \displaystyle C\int_0^1  \frac{e^{-c\frac{4^i}{n(1+w)}}}{\sqrt{w(w+1)}}dw\\
& \leq & \displaystyle Ce^{-c4^i},
\end{array}
\]
where $C,c>0$ only depend on $n$. For $m=0$,
\[
\sum_{0<l<k^2} \frac {e^{-c\frac{4^{i}k^2}{l}}}{l}\leq \int_0^1 e^{-c\frac{4^{i}}u}\frac{du}u\leq Ce^{-c4^{i}}.
\] 
Therefore,
\begin{equation} \label{s1}
\left\vert S_1\right\vert\leq Ce^{-c4^i} \left\Vert f\right\Vert_{L^2(2^{i+1}B\setminus 2^{i}B)}.
\end{equation}
As for $S_2$, Lemma \ref{estim} gives at once 
\begin{equation} \label{s2}
\left\vert S_2\right\vert \leq C e^{-c4^i} \left\Vert f\right\Vert_{L^2(2^{i+1}B\setminus 2^{i}B)},
\end{equation}
where $C,c>0$ only depend on $n$ once more. Finally, for $S_3$, Lemma \ref{estim} provides
\[
\left\vert S_3\right\vert\leq Ck^{2n} \sum_{l>(n+1)k^2} l^{-n-\frac 12} \frac{e^{-c\frac{4^{i}k^2}{l}}}{\sqrt{l}} \left\Vert f\right\Vert_{L^2(2^{i+1}B\setminus 2^{i}B)}.
\]
But one clearly has
\[
\begin{array}{lll}
\displaystyle \sum_{l>(n+1)k^2} l^{-n-\frac 12} \frac{e^{-c\frac{4^{i}k^2}{l}}}{\sqrt{l}} & \leq & \displaystyle \int_{(n+1)k^2}^{+\infty} t^{-n-\frac 12} \frac{e^{-c\frac{4^{i}k^2}{t}}}{\sqrt{t}}dt\\
& = & \displaystyle (4^ik^2)^{-n} \int_{\frac{n+1}{4^i}}^{+\infty} u^{-n}e^{-\frac cu}\frac{du}u\\
& \leq & \displaystyle Ck^{-2n}4^{-in}\int_{0}^{+\infty} u^{-n}e^{-\frac cu}\frac{du}u\leq C4^{-in},
\end{array}
\]
so that, since $k\geq 1$,
\begin{equation} \label{s3}
\left\vert S_3\right\vert\leq C4^{-in}\left\Vert f\right\Vert_{L^2(2^{i+1}B\setminus 2^{i}B)}.
\end{equation}
Summing up the upper estimates (\ref{s1}), (\ref{s2}) and (\ref{s3}) and using (\ref{s1s2s3}), one obtains
\begin{equation} \label{cond1bis3}
\left\Vert T(I-P^{k^2})^nf_i\right\Vert_{L^2(B)}\leq C4^{-in}\left\Vert f\right\Vert_{L^2(2^{i+1}B\setminus 2^{i}B)}.
\end{equation}
The definition of the maximal function and property (\ref{Dbis}) yield
\[
\left\Vert f\right\Vert_{L^2(2^{i+1}B\setminus 2^{i}B)}\leq V^{1/2}(2^{i+1}B)\left({\mathcal M}(\left\vert f\right\vert^2)(x)\right)^{1/2}\leq C2^{(i+1)D/2}V(B)^{1/2} \left({\mathcal M}(\left\vert f\right\vert^2)(x)\right)^{1/2}.
\]
Choosing now $\displaystyle n>\frac D4$ and summing up over $i\geq 1$, one concludes from (\ref{cond1bis1}) and (\ref{cond1bis3}) that
\[
\left\Vert T(I-P^{k^2})^nf\right\Vert_{L^2(B)}\leq C\left(\sum_{i=0}^{+\infty} 2^{i\left(\frac D2-2n\right)}\right) V(B)^{1/2} \left({\mathcal M}(\left\vert f\right\vert^2)(x)\right)^{1/2},
\]
which ends the proof of (\ref{cond1bis}). \hfill\fin\par
\noindent{\bf Proof of (\ref{cond2bis}): }We use the following lemma:
\begin{lem} \label{gradball}
For all $p\in (2,p_0)$, there exists $C,\alpha>0$ such that, for all ball $B\subset\Gamma$ with radius $k$, all integer $i\geq 1$ and all function $f\in L^2(\Gamma)$ supported in $C_i(B)$, and for all $j\in \left\{1,...,n\right\}$ (where $n$ is chosen as above), one has
\[
\left(\frac 1{V(B)^{1/p}}\right)\left\Vert \nabla P^{jk^2}f\right\Vert_{L^p(B)}\leq \frac{Ce^{-\alpha 4^i}}k \frac 1{V(2^{i+1}B)^{1/2}} \left\Vert f\right\Vert_{L^2(\Gamma)}.
\]
\end{lem}
{\bf Proof of Lemma \ref{gradball}: }This proof is very similar to the one of Lemma 3.2 in \cite{acdh}, and we will therefore only indicate the main steps. Consider first the case when $i=1$. If $j=2m$ for some integer $m\geq 0$, (\ref{gradlpoff}) yields
\begin{equation} \label{gradp}
\left\Vert \nabla P^{jk^2}f\right\Vert_{L^p(B)}\leq \frac C{k} \left\Vert P^{mk^2}f\right\Vert_{L^p(\Gamma)}.
\end{equation}
Using (\ref{upper}), and noticing that, by (\ref{D}), for $y\in B$, $V(y,k\sqrt{m})\sim V(B)$, one has, for all $x\in \Gamma$ and all $y\in B$,
\[
p_{mk^2}(x,y)\leq \frac C{V(B)}\exp\left(-c\frac{d^2(x,y)}{mk^2}\right)m(y).
\]
As a consequence, for all $x\in \Gamma$,
\begin{equation} \label{infty}
\left\vert P^{mk^2}f(x)\right\vert\leq \frac C{V^{1/2}(B)}\left\Vert f\right\Vert_{L^2(4B)}.
\end{equation}
The $L^2$ contractivity of $P$ shows that
\begin{equation} \label{contract}
\left\Vert P^{mk^2}f\right\Vert_{L^2(\Gamma)}\leq C\left\Vert f\right\Vert_{L^2(4B)},
\end{equation}
so that, gathering (\ref{infty}) and (\ref{contract}),
\begin{equation} \label{Holder}
\left\Vert P^{mk^2}f\right\Vert_{L^p(\Gamma)}\leq CV(B)^{\frac 1p-\frac 12} \left\Vert f\right\Vert_{L^2(\Gamma)}.
\end{equation}
Finally, (\ref{Holder}) and (\ref{gradp}) yield the conclusion of Lemma \ref{gradball} when $i=1$ and $j=2m$. If $j=2m+1$, argue similarly, writing $j=m+(m+1)$. \par
Consider now the case when $i\geq 2$ and assume that $j=2m$ (one argues similarly if $j=2m+1=m+(m+1)$). Let $\chi_l$ the characteristic function of $C_l(B)$ for all $l\geq 1$. One has, for all $x\in \Gamma$,
\[
\nabla P^{jk^2}f(x)\leq \sum_{l\geq 1}  \nabla P^{mk^2}\chi_lP^{mk^2}f(x)=:\sum_{l\geq 1} g_l(x).
\]
By (\ref{gradlpoff}) and (\ref{Dbis}), 
\[
\begin{array}{lll}
\displaystyle \frac 1{V^{1/p}(B)}\left\Vert g_l\right\Vert_{L^p(B)} & \leq & \displaystyle C \left(\frac{V(2^{l+1}B)}{V(B)}\right)^{1/p} \frac{e^{-c4^l}}{k} \frac 1{V^{1/p}(2^{l+1}B)} \left\Vert P^{mk^2}f\right\Vert_{L^p(2^{l+1}B\setminus 2^lB)}\\
& \leq & \displaystyle C2^{(l+1)D/p} \frac{e^{-c4^l}}{k} \frac 1{V^{1/p}(2^{l+1}B)} \left\Vert P^{mk^2}f\right\Vert_{L^p(2^{l+1}B\setminus 2^lB)}.
\end{array}
\]
Using (\ref{upper}) and arguing as in the proof of Lemma 3.2 in \cite{acdh}, one obtains
\begin{equation} \label{acdh1}
\frac 1{V(2^{l+1}B)}\left\Vert P^{mk^2}f\right\Vert_{L^2(C_l)}^2\leq K_{il}\frac 1{V(2^{i+1}B)} \left\Vert f\right\Vert_{L^2(C_i)}^2
\end{equation}
and, for all $x\in 2^{l+1}B\setminus 2^lB$,
\begin{equation} \label{acdh2}
\left\vert P^{mk^2}f(x)\right\vert\leq K_{il}2^{(i+2)D}\frac 1{V^{1/2}(2^{i+1}B)} \left\Vert f\right\Vert_{L^2(2^{i+1}B\setminus 2^{i}B)},
\end{equation}
where
\[
K_{il}=\left\{
\begin{array}{ll}
Ce^{-c4^{i}} &\mbox{ if }l\leq i-2,\\
C &\mbox{ if }i-1\leq l\leq i+1,\\
Ce^{-c4^l} & \mbox{ if  }l\geq i+2.
\end{array}
\right.
\]
Interpolating between (\ref{acdh1}) and (\ref{acdh2}) therefore yields
\[
\frac 1{V^{1/p}(2^{l+1}B)}\left\Vert P^{mk^2}f\right\Vert_{L^p(C_l)}\leq K_{il}2^{(i+2)D\left(1-\frac 2p\right)} \frac 1{V^{1/2}(2^{i+1}B)} \left\Vert f\right\Vert_{L^2(C_i)}.
\]
Summing up in $l$, one ends the proof of Lemma \ref{gradball} as in \cite{acdh}. \hfill\fin\par
To prove (\ref{cond2bis}), it is enough to show that, if $p\in (2,p_0)$, there exists $C_p>0$ such that, for all $j\in \left\{1,...,n\right\}$, all function $f\in L^2_{loc}(\Gamma)$ with $\nabla f\in L^2_{loc}(\Gamma)$, all ball $B\subset \Gamma$ with radius $k$ and any point $x\in B$,
\[
\frac 1{V^{1/p}(B)}\left\Vert  \nabla P^{jk^2}f \right\Vert_{L^p(B)}\leq C\left({\mathcal M}(\left\vert \nabla f\right\vert^2)\right)^{1/2}(x).
\]
But, since for all $l\geq 0$, $P^l1=1$, one has
\[
\nabla P^lf=\nabla P^l(f-f_{4B}),
\]
so that
\[
\nabla P^{jk^2}f=\sum_{l\geq 1}\nabla P^{jk^2}(\chi_l(f-f_{4B})).
\]
One concludes the proof of (\ref{cond2bis}) as in \cite{acdh}, using the Poincar\'e inequality and Lemma \ref{gradball}. \hfill\fin 
\SE{The Calder\'on-Zygmund decomposition for functions in Sobolev spaces} \label{decompo}
The present section is devoted to the proof of Proposition \ref{CZ}, for which we adapt the proof of Proposition 1.1 in \cite{ac} to the discrete setting. Let  $f\in \dot{E}^{1,p}(\Gamma)$, $\alpha>0$. Consider
$\Omega=\left\lbrace x \in \Gamma : \mathcal{M}(|\nabla f|^{q})(x)>\alpha^{q}\right\rbrace$. If $\Omega=\emptyset$, then set
$$
 g=f\;,\quad b_{i}=0 \, \text{ for all } i\in I
$$
so that (\ref{eg}) is satisfied thanks to the Lebesgue differentiation theorem and the other properties in Proposition \ref{CZ} obviously hold. Otherwise the Hardy-Littlewood maximal theorem gives
\begin{align}
	m(\Omega)&\leq C\alpha^{-p}\|(\nabla f)^{q}\|_{\frac pq}^{\frac pq} \nonumber\\
			& = C \alpha^{-p} \Bigr(\sum_{x} |\nabla f|^{p}(x)m(x)\Bigl) \label{mO}
\\
			&<+\infty. \nonumber
\end{align}
 In particular, $\Omega$ is a proper open subset of $\Gamma$, as $m(\Gamma)=+\infty$ (see Remark \ref{infinite}). Let $(\underline{B_{i}})_{i\in I}$ be a Whitney decomposition of $\Omega$ (\cite{coifman1}). That is, $\Omega$ is the union of the $\underline{B_{i}}$'s, the $\underline{B_{i}}$'s being pairwise disjoint open balls, and there exist two constants $C_{2}>C_{1}>1$, depending only
on the metric, such that, if $F=\Gamma\setminus \Omega$,
\begin{itemize}
\item[1.] the balls $B_{i}=
C_{1}\underline{B_{i}}$ are contained in $\Omega$ and have the bounded overlap property;
\item[2.] for each $i\in I$, $r_{i}=r(B_{i})=\frac{1}{2}d(x_{i},F)$ where $x_{i}$ is 
the center of $B_{i}$;
\item[3.] for each $i\in I$, if $\overline{B_{i}}=C_{2}\underline{B_{i}}$, $\overline{B_i}\cap F\neq \emptyset$ ($C_{2}=4C_{1}$ works).
\end{itemize}
For $x\in \Omega$, denote $I_{x}=\left\lbrace i\in I;\ x\in B_{i}\right\rbrace$. By the bounded overlap property of the balls $B_{i}$, there exists an integer $N$ such that $\sharp I_{x} \leq N$ for all $x\in \Omega$. Fixing $j\in I_{x}$ and using the properties of the $B_{i}$'s, we easily see that $\frac{1}{3}r_{i}\leq r_{j}\leq 3r_{i}$ for all $i\in I_{x}$. In particular, $B_{i}\subset 7B_{j}$ for all $i\in I_{x}$.

Condition (\ref{rb}) is nothing but the bounded overlap property of the $B_{i}$'s  and (\ref{eB}) follows from (\ref{rb}) and  (\ref{mO}). The doubling property and the fact that $\overline{B_{i}} \cap F
\neq \emptyset$ yield:
\begin{equation} \label{f}
\sum_{x\in B_{i}} |\nabla f|^{q}(x)m(x)  \leq
\sum_{x\in \overline{B_{i}}} |\nabla f|^{q}(x)m(x) 
\leq \alpha^{q} V(\overline{B_{i}})
\leq C \alpha^{q}V(B_{i}).
\end{equation}

Let us now define the functions $b_{i}$'s. Let $(\chi_{i})_{i\in I}$ be  a partition of unity of $\Omega$ subordinated to the covering $(\underline{B_{i}})_{i\in I}$, which means that, for all $i\in I$, $\chi_{i}$ is a Lipschitz function supported in $B_{i}$ with
$\displaystyle\|\nabla\chi_{i} \|_{\infty}\leq
\frac{C}{r_{i}}$ and $\displaystyle \sum_{i\in I} \chi_i(x)=1$ for all $x\in \Gamma$ (it is enough to choose $\displaystyle\chi_{i}(x)=
\psi\left(\frac{C_{1}d(x_{i},x)}{r_{i}}\right)\left(\sum_{k}\psi\left(\frac{C_{1}d(x_{k},x)}{r_{k}}\right)\right)^{-1}$, where $\psi\in {\mathcal D}(\R)$, $\psi=1$ on $[0,1]$, $\psi=0$
on $[\frac{1+C_{1}}{2},+\infty)$ and $0\leq \psi\leq 1$). Note that 
$\nabla\chi_{i}$ is supported in $2B_{i}\subset \Omega$.	
We set $b_{i}=(f-f_{B_{i}})\chi_{i}$. It is clear that $\supp b_{i} \subset B_{i}$.
Let us estimate $\sum_{x\in B_{i}} | \nabla b_{i}|^{q}(x)m(x)$. Since
$$
\nabla b_{i}(x)=\nabla ((f-f_{B_{i}})\chi_{i})(x)\leq \max_{y\sim x}\chi_{i}(y) 
\nabla f(x)+|f(x)-f_{B_{i}}|\nabla \chi_{i}(x)
$$
and since $\chi_i(y)\leq 1$ for all $y\in \Gamma$, we get by $(P_{q})$ and (\ref{f}) that
\begin{align*}
\sum_{x\in B_{i}} |\nabla b_{i}|^{q}m(x)
 &\leq
C\left(\sum_{x\in B_{i}}| \nabla f|^{q}(x)m(x)
+\sum_{x\in B_{i}}|f-f_{B_{i}}|^{q}(x)| \nabla \chi_{i}|^{q}(x)m(x) \right)
\\
&\leq C\alpha^{q}V(B_{i})+ C\frac{C^{q}}{r_{i}^{q}} r_{i}^{q}\sum_{x\in 
B_{i}}|\nabla f|^{q}(x)m(x)
\\
&\leq C'\alpha^{q}V(B_{i}).
\end{align*}
Thus (\ref{eb}) is proved. 

Set $ \displaystyle g=f-\sum_{i\in I}b_{i}$. Since the sum is locally finite on $\Omega$,  $g$ is defined everywhere on $\Gamma $ and $g=f$ on $F$. 
 
It remains to prove (\ref{eg}). Since $\displaystyle \sum_{i\in I}\chi_{i}(x)=1$ for all $x\in \Omega$, one has 
$$
g=f\chi_{F}+\sum_{i\in I}f_{B_{i}}\chi_{i}
$$
where $\chi_F$ denotes the characteristic function of $F$. 
%Thus 
%$$ 
%|\nabla g|\leq |\nabla f| \chi_{F}+|\nabla \chi_{F}| |\mathcal{M}f|+|\sum_{i}f_{B_{i}}\nabla \chi_{i}|.
%$$
We will need the following lemma:
\begin{lem} \label{gradg}
There exists $C>0$ such that, for all $j\in I$, all $u\in F\cap 4{B_j}$ and all $v\in B_j$,
\[
\left\vert g(u)-g(v)\right\vert\leq C\alpha d(u,v).
\]
\end{lem}
{\bf Proof: }Since $\displaystyle \sum_{i\in I}\chi_i=1$ on $\Gamma$, one has
\begin{equation} \label{gugv}
\begin{array}{lll}
g(u)-g(v)& = & \displaystyle f(u)-\sum_{i\in I} f_{B_i}\chi_i(v)\\
& = & \displaystyle \sum_{i\in I} \left(f(u)-f_{B_i}\right)\chi_i(v).
\end{array}
\end{equation}
For all $i\in I$ such that $v\in B_i$,
\[
\left\vert f(u)-f_{B_i}\right\vert \leq  \sum_{k=0}^{+\infty} \left\vert f_{B(u,2^{-k}r_i)}-f_{B(u,2^{-k-1}r_i)}\right\vert + \left\vert f_{B(u,r_i)}-f_{B_i}\right\vert.
\]
For all $k\geq 0$, $(P_q)$ yields
\begin{equation} \label{firstterm}
\begin{array}{lll}
\displaystyle \left\vert f_{B(u,2^{-k}r_i)}-f_{B(u,2^{-k-1}r_i)}\right\vert & = & \displaystyle  \frac 1{V(u,2^{-k-1}r_i)}\left\vert \sum_{z\in B(u,2^{-k-1}r_i)} \left(f(z)-f_{B(u,2^{-k}r_i)}\right)m(z)\right\vert \\
& \leq & \displaystyle  \frac C{V(u,2^{-k}r_i)} \sum_{z\in B(u,2^{-k}r_i)} \left\vert f(z)-f_{B(u,2^{-k}r_i)}\right\vert m(z)\\
& \leq & \displaystyle \left(\frac C{V(u,2^{-k}r_i)} \sum_{z\in B(u,2^{-k}r_i)} \left\vert f(z)-f_{B(u,2^{-k}r_i)}\right\vert^q m(z)\right)^{\frac 1q}\\
& \leq & \displaystyle C2^{-k}r_i \left(\frac 1{V(u,2^{-k}r_i)} 
\sum_{z\in B(u,2^{-k}r_i)} \left\vert \nabla f(z)\right\vert^qm(z)\right)^{\frac 1q}\\
& \leq & \displaystyle C2^{-k}r_i \left({\mathcal M}\left(\nabla f\right)^q\right)^{\frac 1q}(u)\\ 
& \leq & \displaystyle C2^{-k}\alpha r_i\leq C2^{-k}\alpha r_j,
\end{array}
\end{equation}
where the penultimate inequality relies on the fact that $u\in F$ and the last one from the fact that $B_i\cap B_j\neq \emptyset$. Moreover, since $u\in 4B_j$,
\[
\begin{array}{lll}
\displaystyle B(u,r_i) & \subset & B(x_j,r_i+d(u,x_j))\\
& \subset & \displaystyle B(x_j,r_i+4r_j)\subset 7B_j.
\end{array}
\]
Since one also has $B_i\subset 7B_j$, one obtains, arguing as before,
\begin{equation} \label{secondterm}
\begin{array}{lll}
\displaystyle \left\vert f_{B(u,r_i)}-f_{B_i}\right\vert & \leq & \displaystyle \left\vert f_{B(u,r_i)}-f_{7B_j}\right\vert + \left\vert f_{7B_j}-f_{B_i}\right\vert\\
& \leq &\displaystyle \frac C{V(7B_j)} \sum_{z\in 7B_j} \left\vert f(z)-f_{7B_j}\right\vert m(z)\\
& \leq & \displaystyle C\alpha r_j.
\end{array}
\end{equation}
It follows from (\ref{firstterm}) and (\ref{secondterm}) that
\[
\left\vert f(u)-f_{B_i}\right\vert \leq C\alpha r_j\leq C\alpha d(u,v),
\]
since
\[
\begin{array}{lll}
r_j & = & \displaystyle \frac 12 d(x_j,F)\leq \frac 12 d(x_j,u)\leq \frac 12 d(x_j,v)+\frac 12 d(v,u)\\
& \leq & \displaystyle \frac 12 r_j+\frac 12 d(v,u).
\end{array}
\]
This ends the proof of Lemma \ref{gradg} because of (\ref{gugv}).\hfill\fin\par
\noindent To prove (\ref{eg}), it is clearly enough to check that $\left\vert g(x)-g(y)\right\vert\leq C\alpha $ for all $x\sim y\in \Gamma$. Let us now prove this fact, distinguishing between three cases:
\begin{itemize}
\item[1.] Assume that $x \in \Omega$. Then, $x\in B_j$ for some $j\in I$, and for all $y\sim x$, $y\in 2B_j\subset \Omega$, so that $\chi_F(x)=\chi_F(y)=0$. It follows that
\[
g(y)-g(x)=\sum_{i\in I} \left(f_{B_i}-f_{B_j}\right)(\chi_i(y)-\chi_i(x)),
\]
so that $\left\vert g(y)-g(x)\right\vert\leq C\sum_{i\in I}\left\vert 
f_{B_{i}}-f_{B_j}\right\vert \nabla\chi_{i}(x):=h(x)$. We claim that 
$|h(x)|\leq C\alpha$. To see this, note that,  for all $i\in I$ such that 
$\nabla\chi_i(x)\neq 0$, we have $|f_{B_{i}}-f_{B_{j}}|\leq Cr_{j} \alpha$.
Indeed, $d(x,B_i)\leq 1$, which easily implies that $r_i\leq 3r_j+1\leq 4r_j$, hence $B_{i} \subset 10B_{j}$. As a consequence, we have, arguing as before again,
\begin{align}
|f_{B_{i}}-f_{10B_{j}}| &\leq
\frac{1}{V(B_{i})}\sum_{y\in B_{i}}|f(y)-f_{10B_{j}}|m(y) \nonumber
\\
&\leq \frac{C}{V(B_{j})}\sum_{y\in 10B_{j}}|f(y)-f_{10B_{j}}|m(y) \nonumber
\\
&\leq Cr_{j}\left(\frac{1}{V(10B_{j})}\sum_{y\in10B_{j}}| \nabla f|^{q}(y)m(y)\right)^{\frac{1}{q}} \nonumber
\\
&\leq Cr_{j}\alpha \label{g}
\end{align}
where we used  H\"{o}lder inequality, $(D)$, $(P_{q})$ and the fact that 
$(|\nabla f|^{q})_{10B_{j}}\leq
\mathcal{M}(|\nabla 
f|)^{q}(z)$ for some $z\in F\cap \overline{B_j}$. Analogously $|f_{10B_{j}}-f_{B_{j}}|\leq Cr_{j}\alpha$. Hence 
\begin{align*}
|h(x)| &=\left\vert \sum_{i\in I;\ x\in 
2B_i}(f_{B_{i}}-f_{B_{j}})\nabla\chi_{i}(x)\right\vert
\\
&\leq C\sum_{i\in I;\ x\in 2B_i}|f_{B_{i}}-f_{B_{j}}|r_{i}^{-1}
\\
&\leq CN\alpha .
\end{align*}
 \item[2.] Assume now that  $x\in F\setminus \partial F$. In this case $\left\vert g(y)-g(x)\right\vert=\left\vert f(y)-f(x)\right\vert\leq C\alpha$ by the definition of $F$.
\item[3.] Assume finally that $x\in \partial F$. 
\begin{itemize}
\item[i.] If $y\in F$, we have $|g(y)-g(x)|=|f(x)-f(y)|\leq C\nabla f(x)\leq C\alpha$.
 \item[ii.] Consider now the case when $y\in \Omega$. There exists $j\in I$ such that $y\in B_j$. Since $x\sim y$, one has $x\in 4B_j$, Lemma \ref{gradg} therefore yields
\[
\left\vert g(x)-g(y)\right\vert \leq C\alpha d(x,y)\leq C\alpha.
\]
%If $x\notin C_2B_{j}$, we choose $z\in  F\cap C_2B_{j}$. Therefore,
%\begin{align*}
%|g(x)-g(y)|& \leq |g(x)-g(z)|+|g(z)-g(y)|
%\\
%&\leq C\alpha d(x,z)+C\alpha d(z,y)
%\\
%&\leq C \alpha d(x,y)+C\alpha r_{j}
%\\
%& \leq C\alpha.
%\end{align*}
%To estimate $|g(x)-g(z)|$ we used the same argument as that of item i. in the proof of Proposition 3.6 in \cite{badr1}, noting that in our discrete case and for the construction of the balls $(B_{i})$, we stop when $r(B_{i})<1$.  For $|g(z)-g(y)|$ we used Lemma \ref{gradg} again. For the last inequality, we use $d(x,y)\leq 1$ and $r_j\leq \frac 12 d(x_j,x)\leq \frac 12 d(x_j,y)+\frac 12 d(y,x)\leq \frac 12 r_j+ \frac 12$.
\end{itemize}
%Then we get $|\nabla g(x)|= \left(\sum_{y\sim x} |g(y)-g(x)|^{2}p(x,y)\right)^{\frac{1}{2}}\leq C\alpha$.
\end{itemize}
Thus the proof of Proposition \ref{CZ} is complete.\hfill\fin
\begin{rem} \label{useful}
It is easy to get the following estimate for the $b_i$'s: for all $i\in I$,
\[
\frac 1{V(B_i)}\left\Vert b_i\right\Vert_1\leq \frac 1{V(B_i)^{1/q}} \left\Vert b_i\right\Vert_q\leq C\alpha r_i.
\]
Indeed, the first inequality follows from H\"older and the fact that $b_i$ is supported in $B_i$. Moreover, by $(P_q)$ and (\ref{f}),
\[
\frac 1{V(B_i)^{1/q}} \left\Vert b_i\right\Vert_q=\frac 
1{V(B_i)^{1/q}}\left\Vert f-f_{B_i}\right\Vert_{L^q(B_i)}\leq  Cr_i\frac 
1{V(B_i)^{1/q}}\left\Vert \nabla f\right\Vert_{L^q(B_i)}\leq C\alpha r_i.
\]
\end{rem}
\SE{An interpolation result for Sobolev spaces} \label{interpo}
To prove Theorem \ref{interpolation}, we will characterize the  $K$ functional of interpolation for homogeneous Sobolev spaces in the following theorem.
 \begin{theo}\label{EKH} Under the same hypotheses as 
Theorem \ref{interpolation} we have that
\begin{itemize}
\item[1.] there exists $C_{1}$ such that for every $f \in \dot{W}^{1,q}(\Gamma)+\dot{W}^{1,\infty}(\Gamma)$ and all $t>0$
$$
 K(f,t^\frac{1}{q},\dot{W}^{1,q},\dot{W}^{1,\infty})\geq 
 C_{1}t^{\frac{1}{q}}\left(|\nabla f|^{q**}\right)^{\frac{1}{q}}(t);
$$
\item[2.] for $ q\leq p<\infty$, there exists $C_{2}$ such that for every $f\in \dot{W}^{1,p}(\Gamma)$  and every $t>0$
$$ 
K(f,t^{\frac{1}{q}},\dot{W}^{1,q},\dot{W}^{1,\infty})\leq C_{2} 
t^{\frac{1}{q}}\left(|\nabla f|^{q**}\right)^{\frac{1}{q}}(t).
$$
\end{itemize}
\end{theo}
{\bf Proof: } We first prove item 1.
Assume that $ f= h+g$ with $ h\in \dot{W}^{1,q},\,g  \in
\dot{W}^{1,\infty}$, we then have 
\begin{align*}
\|h\|_{\dot{W}^{1,q}} + t^{\frac{1}{q}} \|g\|_{\dot{W}^{1,\infty}}
 & \geq \|\nabla h\|_{q} + t^{\frac{1}{q}} \|\nabla g \|_{\infty}
\\
& \geq K(\nabla f,t^{\frac{1}{q}},L^{q},L^{\infty})
\\
&\geq Ct^{\frac{1}{q}}(|\nabla f|^{q**})^{\frac{1}{q}}(t).
\end{align*}
Hence we conclude that $K(f,t^{\frac{1}{q}},\dot{W}^{1,q},\dot{W}^{1,\infty})\geq 
C_{1}t^{\frac{1}{q}}(|\nabla f|^{q**})^{\frac{1}{q}}(t)$.\par
\noindent We prove now item 2. Let $f\in \dot{W}^{1,p},\, q\leq p<\infty$.
Let $t>0$, we consider the Calder\'{o}n-Zygmund decomposition of $f$ given by Proposition \ref{CZ} with
$\alpha=\alpha(t)=\Bigl(\mathcal{M}(|\nabla f|)^{q}\Bigr)^{*\frac{1}{q}}(t)$. Thus we have $ \displaystyle f=\sum_{i\in I}b_{i}+g=b+g $ where
$(b_{i})_{i\in I},\,g$ satisfy the properties of the proposition. We have the estimate
\begin{align*}
\| \nabla b \|_{q}^{q}&\leq \sum_{x\in \Gamma}\left(\sum_{i\in I}
|\nabla b_{i}|\right)^{q}(x)m(x)
\\
&\leq CN \sum_{i\in I}\sum_{x\in B_{i}}
|\nabla b_{i}|^{q}(x)m(x)
\\
&\leq
C\alpha^{q}(t)\sum_{i\in I}V(B_{i})
\\
&\leq C\alpha^{q}(t)m(\Omega),
\end{align*}
where the $B_i$'s are given by Proposition \ref{CZ} and $\Omega$ is defined as in the proof of Proposition \ref{CZ}. The last inequality follows from the fact that
$\displaystyle \sum_{i\in I}\chi_{B_{i}}\leq N$ and
$\Omega=\underset{i}{\bigcup}B_{i}$. Hence $\|\nabla b\|_{q}\leq
C\alpha(t)m(\Omega)^{\frac{1}{q}}$. Moreover, since $(\mathcal{M}f)^{*}\sim f^{**}$ (see \cite{bs}, Chapter 3, Theorem 3.8), we obtain  
\[
\alpha(t)=\left(\mathcal{M}(|\nabla f|)^{q}\right)^{*\frac{1}{q}}(t)\leq 
C\left(|\nabla f|^{q**}\right)^{\frac{1}{q}}(t).
\]
Hence, also noting that $m(\Omega)\leq t$ (see \cite{bs}, Chapter 2, Proposition 1.7), we get $K(f,t^{\frac{1}{q}},\dot{W}^{1,q},\dot{W}^{1,\infty})\leq
Ct^{\frac{1}{q}}|\nabla f|^{q**\frac{1}{q}}(t)$ for all $t>0$ and obtain the desired inequality. \hfill\fin\par
{\bf Proof of Theorem \ref{interpolation}: } The proof follows directly from Theorem \ref{EKH}. Indeed, item 1. of Theorem \ref{EKH} gives us that  $(\dot{W}^{1,q},\dot{W}^{1,\infty})_{1-\frac{q}{p},p} \subset \dot{W}^{1,p}$ and 
$\|f\|_{\dot{W}^{1,p}}\leq C\|f\|_{1-\frac{q}{p},p}$, while  item 2. gives us that $\dot{W}^{1,p}\subset (\dot{W}^{1,q},\dot{W}^{1,\infty})_{1-\frac{q}{p},p}$ and $\|f\|_{1-\frac{q}{p},p}\leq C\|f\|_{\dot{W}^{1,p}}$. Hence 
$\dot{W}^{1,p}= (\dot{W}^{1,q},\dot{W}^{1,\infty})_{1-\frac{q}{p},p}$ with equivalent norms. \hfill\fin
\SE{The proof of $(RR_p)$ for $p<2$} \label{reverseriesz}
In view of Theorem \ref{interpolation} and since $(RR_2)$ holds, it is enough, for the proof of 
Theorem \ref{reverserieszpleq2}, to establish (\ref{enoughestimate}).

{\bf Proof of (\ref{enoughestimate}): }We follow the proof of (1.9) in \cite{ac}. Consider such an $f$ and fix $\lambda>0$. Perform the 
Calderon-Zygmund decomposition of $f$ given by Proposition \ref{CZ}. We also 
use the following expansion of $(I-P)^{1/2}$:
\begin{equation} \label{expansionbis}
(I-P)^{1/2}=\sum_{k=0}^{+\infty} a_{k}(I-P)P^k
\end{equation}
where the $(a_{k})$'s were already considered in Section \ref{riesz}. For each $i\in I$, pick the integer $k\in \Z$ such that 
$2^k\leq r(B_{i})<2^{k+1}$ and define $r_{i}=2^k$. We split the 
expansion (\ref{expansionbis}) into two parts:
\[
(I-P)^{1/2}=\sum_{k=0}^{r_{i}^{2}} a_{k}(I-P)P^k+ 
\sum_{k=r_{i}^{2}+1}^{+\infty} a_{k}(I-P)P^k:=T_{i}+U_{i}.
\]
We first claim that
\begin{equation} \label{claim1}
m\left(\left\{x\in \Gamma;\ \left\vert 
(I-P)^{1/2}g(x)\right\vert>\lambda\right\}\right)\leq \frac 
C{\lambda^q}\left\Vert \nabla f\right\Vert_{q}^q.
\end{equation}
Indeed, one has
\[
\begin{array}{lll}
\displaystyle m\left(\left\{x\in \Gamma;\ \left\vert 
(I-P)^{1/2}g(x)\right\vert>\lambda\right\}\right) & \leq & 
\displaystyle \frac C{\lambda^2}\left\Vert 
(I-P)^{1/2}g\right\Vert_{2}^2\\
& = & \displaystyle \frac C{\lambda^2}\left\Vert 
\nabla g\right\Vert_{2}^2,
\end{array}
\]
and since $\nabla g\leq C\lambda$ on $\Gamma$ and $\left\Vert 
\nabla g\right\Vert_{q}\leq C\left\Vert \nabla f\right\Vert_{q}$, we obtain
\[
\left\Vert 
\nabla g\right\Vert_{2}^2 \leq C\lambda^{2-q} \left\Vert 
\nabla g\right\Vert_{q}^q\leq C\lambda^{2-q}\left\Vert 
\nabla f\right\Vert_{q}^q,
\]
which ends the proof of (\ref{claim1}). \par
We now claim that, for some constant $C>0$,
\begin{equation} \label{claim2}
m\left(\left\{x\in \Gamma;\ \left\vert \sum_{i\in 
I}T_{i}b_{i}(x)\right\vert>\lambda\right\}\right)\leq \frac 
C{\lambda^q}\left\Vert \nabla f\right\Vert_{q}^q.
\end{equation}    
To prove (\ref{claim2}), write
\begin{equation} \label{claim2step1}
m\left(\left\{x\in \Gamma;\ \left\vert \sum_{i\in I}T_{i}b_{i}(x)\right\vert>\lambda\right\}\right)\leq m\left(\bigcup_{i} 4B_i\right)+ m\left(\left\{x\notin \bigcup_{i} 4B_i;\ \left\vert \sum_{i\in I}T_{i}b_{i}(x)\right\vert>\lambda\right\}\right).
\end{equation}
Observe first that, by (\ref{D}) and Proposition \ref{CZ},
\[
m\left(\bigcup_{i} 4B_i\right) \leq C\sum_{i\in I} V(4B_i)\leq  \frac 
C{\lambda^q}\left\Vert \nabla f\right\Vert_{q}^q.
\]
As far as the second term in the right-hand side of (\ref{claim2step1}) is concerned, it can be estimated by
\[
m\left(\left\{x\notin \bigcup_{i} 4B_i;\ \left\vert \sum_{i\in I}T_{i}b_{i}(x)\right\vert>\lambda\right\}\right)\leq \frac 1{\lambda^2} \sum_{x\in \Gamma} \left\vert  \sum_{i\in I}\chi_{\Gamma\setminus 4B_i}(x)T_{i}b_{i}(x)\right\vert^2m(x).
\]
Arguing as in \cite{ac, bk, hm}, we estimate this last quantity by duality. Fix a function $u\in L^2(\Gamma,m)$ with $\left\Vert u\right\Vert_2=1$. One has
\[
\left\vert \sum_{x\in \Gamma} \sum_{i\in I}\chi_{\Gamma\setminus 4B_i}(x)T_{i}b_{i}(x) u(x)m(x)\right\vert \leq \sum_{i\in I} \sum_{j=2}^{+\infty} A_{i,j}
\]
where, for all $i\in I$ and all $j\geq 2$,
\[
A_{i,j}:=\sum_{x\in 2^{j+1}B_i\setminus 2^jB_i} \left\vert T_ib_i(x)\right\vert \left\vert u(x)\right\vert m(x).
\]
If $i,j$ are fixed, since $(I-P)b_i$ is supported in $2B_i$,
\[
\begin{array}{lll}
\displaystyle \left\Vert T_ib_i\right\Vert_{L^2(2^{j+1}B_i\setminus 2^jB_i)} & \leq & \displaystyle \sum_{k=0}^{r_{i}^{2}} \left\vert a_{k}\right\vert \left\Vert (I-P)P^kb_i\right\Vert_{L^2(2^{j+1}B_i)\setminus 2^jB_i)}\\
& = & \displaystyle \sum_{k=1}^{r_{i}^{2}} \left\vert a_{k}\right\vert \left\Vert (I-P)P^kb_i\right\Vert_{L^2(2^{j+1}B_i)\setminus 2^jB_i)}
\end{array}
\]
Given $1\leq k\leq r_i^2$, one has, for all $x\in 2^{j+1}B_i\setminus 2^jB_i$, using (\ref{timederivative}),
\[
\left\vert (I-P)P^kb_i(x)\right\vert\leq \sum_{y\in B_i} \left\vert p_k(x,y)-p_{k+1}(x,y)\right\vert \left\vert b_i(y)\right\vert \leq \sum_{y\in B_i} \frac C{kV(y,\sqrt{k})} e^{-c\frac{d^2(x,y)}k} \left\vert b_i(y)\right\vert m(y).
\]
Using (\ref{Dbis}) and arguing exactly as in \cite{ac} (relying, in particular, on Remark \ref{useful}), we obtain
\[
\left\Vert (I-P)P^kb_i\right\Vert_{L^2(2^{j+1}B_i\setminus 2^jB_i)} \leq C\frac{r_i}k \left(\frac{r_i}{\sqrt{k}}\right)^{2D} e^{-c\frac{4^jr_i^2}k} V^{1/2} (2^{j+1}B_i)\lambda.
\]
Since 
\[
a_k\sim \frac 1{\sqrt{k\pi}}
\]
(see Appendix), it follows that
\[
\left\Vert T_ib_i\right\Vert_{L^2(2^{j+1}B_i\setminus 2^jB_i)}  \leq Ce^{-c4^j}V^{1/2} (2^{j+1}B_i)\lambda.
\]
One concludes, as in \cite{ac}, that (\ref{claim2}) holds.\par
What remains to be proved is that
\begin{equation} \label{claim3}
m\left(\left\{x\in \Gamma;\ \left\vert \sum_{i\in 
I}U_{i}b_{i}(x)\right\vert>\lambda\right\}\right)\leq \frac 
C{\lambda^q}\left\Vert \nabla f\right\Vert_{q}^q.
\end{equation}    
Define, for all $j\in \Z$,
\[
\beta_j=\sum_{i\in I;\ r_i=2^j} \frac{b_i}{r_i},
\]
so that, for all $j\in \Z$,
\[
\sum_{i\in I;\ r_i=2^j} b_i=2^j \beta_j.
\]
One has
\[
\begin{array}{lll}
\displaystyle \sum_{i\in I} U_ib_i &= & \displaystyle \sum_{i\in I} \sum_{k>r_i^2} a_k(I-P)P^kb_i\\
& = & \displaystyle \sum_{k>0} a_k(I-P)P^k\sum_{i\in I;\ r_i^2<k} b_i\\
& =& \displaystyle \sum_{k>0} a_k(I-P)P^k\sum_{i\in I;\ r_i^2=2^{2j}<k} b_i \\
& = & \displaystyle \sum_{k>0} a_k(I-P)P^k\sum_{j;\ 4^j<k} 2^j\beta_j.
\end{array}
\]
For all $k>0$, define
\[
f_k=\sum_{j;\ 4^j<k} \frac{2^j}{\sqrt{k}}\beta_j.
\]
It follows from the previous computation and Theorem \ref{LP} that
\[
\left\Vert \sum_{i\in I} U_ib_i\right\Vert_q\leq C\left\Vert \left(\sum_{k=1}^{+\infty} \frac 1k \left\vert f_k\right\vert^2\right)^{1/2}\right\Vert_q.
\]
To see this, we estimate the left-hand side of this inequality by duality, as in \cite{ac} and use the fact that $\displaystyle \left\vert a_k\right\vert\leq \frac C{\sqrt{k}}$ for all $k\geq 1$. Since, by Cauchy-Schwarz,
\[
\left\vert f_k\right\vert^2\leq 2\sum_{j;\ 4^j<k} \frac{2^j}{\sqrt{k}} \left\vert \beta_j\right\vert^2,
\]
one obtains
\[
\left\Vert \left(\sum_{k=1}^{+\infty} \frac 1k \left\vert f_k\right\vert^2\right)^{1/2}\right\Vert_q\leq \left\Vert \left(\sum_{k\in \Z} \left\vert \beta_k\right\vert^2\right)^{1/2}\right\Vert_q.
\]
By the bounded overlap property,
\[
\left\Vert \left(\sum_{k\in \Z} \left\vert \beta_k\right\vert^2\right)^{1/2}\right\Vert_q^q\leq C\sum_{x\in \Gamma} \sum_{i\in I} \frac{\left\vert b_i(x)\right\vert^q}{r_i^q} m(x),
\]
so that, using Remark \ref{useful}, one obtains
\[
\sum_{x\in \Gamma} \sum_{i\in I} \frac{\left\vert b_i(x)\right\vert^q}{r_i^q} m(x)\leq C\lambda^q\sum_{i\in I} V(B_i).
\]
As a conclusion,
\[
m\left(\left\{x\in \Gamma;\ \left\vert \sum_{i\in 
I}U_{i}b_{i}(x)\right\vert>\lambda\right\}\right)\leq C\sum_{i\in I} 
V(B_i)\leq \frac C{\lambda^q} \left\Vert \nabla f\right\Vert_q^q,
\]
which is exactly (\ref{claim3}). The proof of (\ref{enoughestimate}) is therefore complete. \hfill\fin
\section{Riesz transforms and harmonic functions} \label{rieszharmonic}
Let us now prove Theorem \ref{rpharmonic}. The proof goes through a property analogous to $(\Pi_p)$ in \cite{ac}, the statement of which requires a notion of discrete differential.
\subsection{The discrete differential and its adjoint}
To begin with, for any $\gamma=(x,y),\gamma^{\prime}=(x^{\prime},y^{\prime})\in E$ (recall that $E$ denotes the set of edges in $\Gamma$), set
\[
d(\gamma,\gamma^{\prime})=\max(d(x,x^{\prime}),d(y,y^{\prime})).
\] 
It is straightforward to check that $d$ is a distance on $E$. We also define a measure on subsets on $E$. For any $A\subset E$, set
\[
\mu(A)=\sum_{(x,y)\in A} \mu_{xy}.
\]
We claim that $E$, equipped with the metric $d$ and the measure $\mu$, is a space of homogeneous type. Indeed, let $\gamma=(a,b)\in E$ and $r>0$. Assume first that $r\geq 5$. Then, by (\ref{D}),
\[
\mu(B(\gamma,2r))=\sum_{d(x,a)<2r,\ d(y,b)<2r} \mu_{xy}\leq \sum_{d(x,a)<2r} \sum_{y\in \Gamma} \mu_{xy}=V(a,2r)\leq CV\left(a,\frac r{100}\right).
\]
But
\[
V\left(a,\frac r{100}\right)=\sum_{d(x,a)<\frac r{100}} \sum_{d(y,x)\leq 1} \mu_{xy}\leq \sum_{d(x,a)<\frac r2}\sum_{d(y,b)<\frac r2} \mu_{xy}=\mu\left(B\left(\gamma,\frac r2\right)\right),
\]
since, when $d(x,a)<\frac r{100}$ and $d(y,x)\leq 1$, then $d(y,b)<2+\frac r{100}\leq \frac r2$. \par
\noindent Assume now that $r<5$. One has, using (\ref{D}) again,
\[
\mu(B(\gamma,2r))\leq V(a,2r)\leq V(a,10)\leq CV\left(a,\frac 12\right)=Cm(a)\leq C^{\prime}\mu_{ab}\leq C^{\prime} \mu(B(\gamma,r)),
\]
since, whenever $x\sim y$, one has $\alpha m(x)\leq \mu_{xy}$ by (\ref{deltaalpha}). The claim is therefore proved.\par
\noindent We can then define $L^p$ spaces on $E$ in the following way. For $1\leq p< +\infty$, say that a function $F$ on $E$ belongs to $L^p(E)$ if and only if $F$ is antisymmetric (which means that $F(x,y)=-F(y,x)$ for all $(x,y)\in E$) and
\[
\left\Vert F\right\Vert_{L^p(E)}^p:=\frac 12\sum_{(x,y)\in E} \left\vert F(x,y)\right\vert^p\mu_{xy}<+\infty.
\]
Observe that the $L^2(E)$-norm derives from the scalar product 
\[
\langle F,G\rangle_{L^2(E)}=\frac 12\sum_{x,y\in \Gamma} F(x,y)G(x,y)\mu_{xy}.
\]
Finally, say that $F\in L^{\infty}(E)$ if and only if $F$ is antisymmetric and 
\[
\left\Vert F\right\Vert_{L^{\infty}(E)}:=\frac 12\sup_{(x,y)\in E} \left\vert F(x,y)\right\vert<+\infty.
\]
Our notion of discrete differential is the following one: for any function $f$ on $\Gamma$ and any $\gamma=(x,y)\in E$, define
\[
df(\gamma)=f(y)-f(x).
\]
The function $df$ is clearly antisymmetric on $E$ and is related to the length of the gradient of $f$. More precisely, it is not hard to check that, if (\ref{deltaalpha}) holds, then for all $p\in [1,+\infty]$ and all function $f$ on $\Gamma$,
\[
\left\Vert df\right\Vert_{L^p(E)}\sim \left\Vert \nabla f\right\Vert_{L^p(\Gamma)}.
\]
Indeed, if $1\leq p<+\infty$, for all function $f$ and all $x\in \Gamma$,
\[
\begin{array}{lll}
\displaystyle \left\vert \nabla f(x)\right\vert^p &\sim & \displaystyle \left(\sum_{y\sim x} p(x,y)\left\vert f(y)-f(x)\right\vert\right)^p\\
& \sim & \displaystyle \sum_{y\sim x} p^p(x,y)\left\vert f(y)-f(x)\right\vert^p\\
& \sim & \displaystyle \sum_{y\sim x} p(x,y)\left\vert f(y)-f(x)\right\vert^p
\end{array}
\]
where the last line is due to (\ref{deltaalpha}). As a consequence,
\[
\begin{array}{lll}
\displaystyle \left\Vert \nabla f\right\Vert_{L^p(\Gamma)}^p & \sim & \displaystyle \sum_{x\in \Gamma} \sum_{y\sim x} p(x,y)\left\vert f(y)-f(x)\right\vert^pm(x)\\
& \sim & \displaystyle \sum_{x,y\in \Gamma} \left\vert df(x,y)\right\vert^p\mu_{xy}\\
& = & \displaystyle \left\Vert df\right\Vert_{L^p(E)}^p.
\end{array}
\]
The case when $p=+\infty$ is analogous and even easier. We could therefore reformulate properties (\ref{rp}) and (\ref{rrp}) replacing $\left\Vert \nabla f\right\Vert_{L^p(\Gamma)}$ by $\left\Vert df\right\Vert_{L^p(E)}$.\par
\noindent Besides $d$, we also consider its adjoint in $L^2$. If $df\in L^2(E)$ and $G$ is any (antisymmetric) function in $L^2(E)$ such that the function $x\mapsto \sum_{y} p(x,y)G(x,y)$ belongs to $L^2(\Gamma)$, one has
\[
\begin{array}{lll}
\displaystyle \langle df,G\rangle_{L^2(E)} & = & \displaystyle \frac 12\sum_{x,y\in \Gamma}  df(x,y)G(x,y)\mu_{xy}\\
& = & \displaystyle \frac 12\sum_{x,y\in \Gamma}  f(y)G(x,y)\mu_{xy} - \frac 12\sum_{x,y\in \Gamma}  f(x)G(x,y)\mu_{xy}\\
& =& \displaystyle -\sum_{x,y\in \Gamma}  f(x)G(x,y)\mu_{xy}\\
& = & \displaystyle -\sum_{x\in \Gamma} f(x)\left(\sum_{y\in \Gamma} p(x,y)G(x,y)\right)m(x).
\end{array}
\] 
Thus, if we define
\[
\delta G(x)=\sum_{y} p(x,y)G(x,y)
\]
for all $x\in \Gamma$, it follows that 
\[
\langle df,G\rangle_{L^2(E)}=-\langle f,\delta G\rangle_{L^2(\Gamma)}
\]
whenever $f\in L^{2}(\Gamma)$, $dF\in L^2(E)$, $G\in L^2(E)$ and $\delta G\in L^2(\Gamma)$.  Notice also that $I-P=-\delta d$.\par
\medskip
\noindent The following lemma, very similar to Lemma 4.2 in \cite{auschcou}, holds:
\begin{lem} \label{inhomogeneous}
Assume that (\ref{D}), (\ref{deltaalpha}) and (\ref{diagupper}) hold. There exists $C>0$ such that, for all ball $B$ and all function $f\in L^2(\Gamma)$ supported in $B$, there exists a unique function $h\in W^{1,2}_0(B)$ such that 
\begin{equation} \label{sol}
(I-P)h=f\mbox{ in }\Gamma
\end{equation}
and $h$ satisfies
\[
\left\Vert h\right\Vert_{W^{1,2}(\Gamma)}\leq C\left\Vert f\right\Vert_{L^2(\Gamma)}.
\]
\end{lem}
{\bf Proof: } This proof relies on a Sobolev inequality, which will be used again in the proof of Theorem \ref{rpharmonic} and reads as follows: there exist $\nu\in (0,1)$ and $C>0$ such that, for all ball $B$ with radius $r>\frac 12$ and all function $f$ supported in $B$, 
\begin{equation} \label{sob}
\left\Vert f\right\Vert_{q}\leq CrV(B)^{-\frac{\nu}2}\left\Vert \nabla f\right\Vert_2
\end{equation}
with $q=\frac 2{1-\nu}$. This inequality is actually equivalent to a relative Faber-Krahn inequality, which is itself equivalent to the conjunction of (\ref{D}) and (\ref{diagupper}), see \cite{cougri, grigo, carron, coulhonlp, bcls, delmoelliptic}. \par
\noindent Let $B$ and $f$ as in the statement of Lemma \ref{inhomogeneous}. Since $I-P=-\delta d$, (\ref{sol}) is equivalent to
\[
\langle dh,dv\rangle_{L^2(E)}=\langle f,v\rangle_{L^2(\Gamma)}
\]
for all $v\in W^{1,2}_0(B)$. For all $u,v\in W^{1,2}_0(B)$, set ${\mathcal B}(u,v)=\langle du,dv\rangle_{L^2(E)}$. It is obvious that ${\mathcal B}$ is a continuous bilinear form on $W^{1,2}_0(B)$. Moreover, for all $u\in W^{1,2}_0(B)$,
\[
{\mathcal B}(u,u)=\left\Vert du\right\Vert_{L^2(E)}^2\geq c\left\Vert u\right\Vert_{W^{1,2}_0(B)}^2,
\]
by (\ref{sob}) (see also Lemma 4.1 in \cite{auschcou}). The conclusion of Lemma \ref{inhomogeneous} follows then from the Lax-Milgram theorem.\hfill\fin\par

\medskip

\noindent Let $F\in L^2(E)$. It is easy to check that $\delta F\in L^2(\Gamma)$ and
\begin{equation} \label{easytocheck}
\left\Vert \delta F\right\Vert_{L^2(\Gamma)}\leq \left\Vert F\right\Vert_{L^2(E)}.
\end{equation}
Indeed, for all $g\in L^2(\Gamma)$,
\[
\begin{array}{lll}
\displaystyle \left\vert \langle \delta F,g\rangle_{L^2(\Gamma)} \right\vert &= & \displaystyle \left\vert \sum_{x,y\in \Gamma}p(x,y)F(x,y)g(x)m(x)\right\vert\\
& = & \displaystyle \left\vert \sum_{x,y\in \Gamma} F(x,y)g(x)\mu_{xy}\right\vert \\
& \leq & \displaystyle \left(\sum_{x,y\in \Gamma} \left\vert F(x,y)\right\vert^2\mu_{x,y}\right)^{1/2} \left(\sum_{x\in \Gamma} \left\vert g(x)\right\vert^2 m(x)\right)^{1/2}.
\end{array}
\]
As a consequence of Lemma \ref{inhomogeneous}, for all $F\in L^2(E)$ with bounded support, there exists a unique function $f\in W^{1,2}(\Gamma)$ such that $(I-P)f=\delta F$. Since functions in $L^2(E)$ with bounded support are dense in $L^2(E)$, we can therefore extend the operator $d(I-P)^{-1}\delta$ to an $L^2(E)$-bounded operator.
\subsection{The proof of Theorem \ref{rpharmonic}}
\noindent For all $1\leq p<+\infty$, say that $(\Pi_p)$ holds if and only if there exists $C_p>0$ such that, for all $F\in L^p(E)\cap L^2(E)$,
\begin{equation} \label{pip} \tag{$\Pi_p$}
\left\Vert d(I-P)^{-1}\delta F\right\Vert_{L^p(E)}\leq C_p\left\Vert F\right\Vert_{L^p(E)}.
\end{equation}
Since $L^2(E)\cap L^p(E)$ is dense in $L^p(E)$, if $(\Pi_p)$ holds, the operator $d(I-P)^{-1}\delta$ extends to a bounded operator from $L^p(E)$ to itself.\par
\noindent Let us now turn to the proof of Theorem \ref{rpharmonic}. Let $p_0>2$ and $q\in (2,p_0)$. Denote by $(2^{\prime})$ the following property: 
\begin{equation} \label{2prime} \tag{$2^{\prime}$}
\mbox{for all }p\in (2,q), (\Pi_p) \mbox{ holds}.
\end{equation}
We show that, for some $p_0>2$, if $q\in (2,p_0)$, then $2.\Rightarrow 2^{\prime}.\Rightarrow 1.\Rightarrow 2.$\par

\medskip

\noindent{\bf Proof of $2.\Rightarrow 2^{\prime}.$ }In order to apply Theorem 2.3 in \cite{ac}, observe first that $E$, equipped with the metric $d$ and the measure $\mu$, is a space of homogeneous type. Let $2<p<\widetilde{p}<q$. Consider $F\in 
L^2(E)\cap L^p(E)$ with bounded support included in $E\setminus 64B$ where $B$ is a ball in 
$E$ centered at $\gamma=(a,b)$ and with radius $r$. Lemma \ref{inhomogeneous} and (\ref{easytocheck}) therefore yield a function $h\in W^{1,2}(\Gamma)$ such that $(I-P)h=\delta F$ with 
$\left\Vert h\right\Vert_{W^{1,2}(\Gamma)}\leq C\left\Vert \delta F\right\Vert_{L^2(\Gamma)}\leq C\left\Vert F\right\Vert_{L^2(E)}$. \par
\noindent If $r\geq \frac 1{16}$, then the function $h$ is harmonic in $B(a,32r)$. Indeed, if $x\in B(a,32r)\setminus \partial B(a,32r)$, 
\[
(I-P)h(x)=\delta F(x)=\sum_{y\sim x} p(x,y)F(x,y).
\]
When $x\in B(a,32r)$ and $y\sim x$, $d(y,b)\leq d(x,a)+2\leq 64r$, so that 
$F(x,y)=0$.  It follows from $\left(RH_{\widetilde{p}}\right)$ that
\[
\left(\frac 1{V(B)} \sum_{x\in B} \left\vert \nabla 
h(x)\right\vert^{\widetilde{p}}m(x)\right)^{\frac 1{\widetilde{p}}} \leq 
C\left(\frac 1{V(16B)} \sum_{x\in 16B} \left\vert \nabla h(x)\right\vert^{2}m(x)\right)^{\frac 12}.
\]
If $r<\frac 1{16}$, $B=16B$ and the same inequality holds. This shows that 
the operator $T$ defined by $TF=\nabla (I-P)^{-1}\delta F$ 
for all $F$ with bounded support in $E$, clearly satisfies the assumptions of Theorem 2.3 in \cite{ac}, and this theorem therefore yields
\begin{equation} \label{Tbounded}
\left\Vert TF\right\Vert_{L^p(E)}\leq C_p\left\Vert F\right\Vert_{L^p(E)}
\end{equation}
for all $F$ with bounded support in $E$. Since the space of antisymmetric functions on $E$ with bounded support 
is dense in $L^p(E)$, (\ref{Tbounded}) holds for all $F\in L^p(E)$, which exactly means that $(\Pi_p)$ holds.\hfill\fin\par

\medskip

\noindent{\bf Proof of $2^{\prime}. \Rightarrow 1.$ } By Theorem \ref{reverserieszpleq2} and Proposition \ref{selfimprove}, there exists $\varepsilon>0$ such that $(RR_q)$ holds for all $q\in (2-\varepsilon,2)$. It is therefore enough to check that the conjunction of $(\Pi_p)$ and $(RR_{p^{\prime}})$ implies $(R_p)$, with $\frac 1p+\frac 1{p^{\prime}}=1$. But, if $f\in L^p(\Gamma)\cap L^2(\Gamma)$ and $G\in L^{p^{\prime}}(E)\cap L^2(E)$,
\[
\begin{array}{lll}
\displaystyle \left\vert \langle d(I-P)^{-1/2}f,G\rangle_{L^2(E)}\rangle\right\vert & = & \displaystyle \left\vert \langle (I-P)^{-1/2}f,\delta G\rangle_{L^2(\Gamma)}\right\vert\\
& = & \displaystyle \left\vert \langle f,(I-P)^{-1/2}\delta G\rangle_{L^2(\Gamma)}\right\vert\\
& \leq & \displaystyle \left\Vert f\right\Vert_{L^p(\Gamma)} \left\Vert (I-P)^{-1/2}\delta G\right\Vert_{L^{p^{\prime}}(\Gamma)}\\
& = & \displaystyle   \left\Vert f\right\Vert_{L^p(\Gamma)} \left\Vert (I-P)^{1/2}(I-P)^{-1}\delta G\right\Vert_{L^{p^{\prime}}(\Gamma)}\\
& \leq & \displaystyle   \left\Vert f\right\Vert_{L^p(\Gamma)} \left\Vert d(I-P)^{-1}\delta G\right\Vert_{L^{p^{\prime}}(\Gamma)}\\
& \leq & \displaystyle C\left\Vert f\right\Vert_{L^p(\Gamma)} \left\Vert G\right\Vert_{L^{p^{\prime}}(\Gamma)},
\end{array}
\]
which ends the proof. \hfill\fin\par

\medskip

\noindent{\bf Proof of $1.\Rightarrow 3.$ } Assume now that $(R_p)$ holds for some $p\in (2,q)$. 
%and that there exist $c>0$ and $d>2$ such that
%\begin{equation} \label{minvol}
%V(x,r)\geq cr^d.
%\end{equation}
%Under assumption (\ref{minvol}), (\ref{diagupper}) yields
%\begin{equation} \label{diagupperbis}
%p_k(x,x)\leq \frac Ck^{-d/2}
%\]
%for all $x\in \Gamma$. Then (\ref{diagupper}) implies (and is actually equivalent to) the following Sobolev inequality:
%\begin{equation} \label{sob}
%\left\Vert f\right\Vert_{L^{\frac{2d}{d-2}}(\Gamma)}\leq C\left\Vert \nabla f\right\Vert_{L^2(\Gamma)}
%\end{equation}
%whenever $f\in \cdot{W}^{1,2}(\Gamma)$ (see \cite{carron, cks, coulhonsurvey, cl, cs, grigo, varo}. As a consequence of (\ref{sob}), there exists $C>0$ such that, for all $h\in L^2(\Gamma)$,
%\begin{equation} \label{rieszpotential}
%\left\Vert (I-P)^{-1/2}h\right\Vert_{\frac{2d}{d-2}}\leq C\left\Vert h\right\Vert_{2}.
%\end{equation}
%Indeed, $\nabla (I-P)^{-1/2}h\right\Vert_2\leq C\left\Vert h\right\Vert_2$.\par
Let $B$ be a ball with center $x_0$ and radius $k$ and $u$ a 
function harmonic in $32B$, and fix a function $\varphi$ supported in $3B$, 
equal to $1$ in $2B$ and satisfying $0\leq \varphi\leq 1$ and $\left\Vert 
\nabla \varphi\right\Vert_{\infty}\leq \frac Ck$. Up to an additive constant, one may assume that the mean value of $u$ in $16B$ is $0$. In order to  control 
the left-hand side of (\ref{reverseholder}), it suffices to estimate 
$\sum_{x\in B} \left\vert \nabla (u\varphi)(x)\right\vert^p m(x)$. \par
%\noindent To that purpose, we follow the appraoch of Shen in \cite{shen}, Lemmata 2.4 and 2.7. A calculation shows that, for all $x\in \Gamma$,
%\begin{equation} \label{deltaformula}
%\begin{array}{lll}
%\displaystyle(I-P)(u\varphi)(x) &= & \displaystyle \sum_{y\in \Gamma} p(x,y)(u\varphi(x)-u\varphi(y))\\
%& = & \displaystyle \sum_{y\in \Gamma} p(x,y)u(x)(\varphi(x)-\varphi(y))+\sum_{y\in \Gamma} p(x,y) (u(x)-u(y))(\varphi(y)-\varphi(x))\\
%& + & \displaystyle \sum_{y\in \Gamma} p(x,y)(u(x)-u(y))\varphi(x)\\
%& := & \displaystyle v_1(x)+v_2(x)+v_3(x).
%\end{array}
%\end{equation}
%For all $x\in \Gamma$, $v_3(x)=0$ since $u$ is harmonic in $3B$. Since $(I-P)(u\varphi)=v_1+v_2$, $u\varphi\in W^{1,2}(\Gamma)$ and $v_1,v_2\in L^2(\Gamma)$, Lemma \ref{inhomogeneous} shows that $u\varphi=(I-P)^{-1}(v_1+v_2)$. Assumption $(R_p)$ therefore yields
%\[
%\left\Vert \nabla(u\varphi)\right\Vert_p\leq C_p \left\Vert (I-P)^{-1/2}(v_1+v_2)\right\Vert_p.
%\]
%Let $p\in (2,p_0)$. Thanks to assumption $(RR_p)$ and (\ref{rieszpotential}), we therefore obtain
%\[
%\left\Vert \nabla(u\varphi)\right\Vert_p\leq C \left\Vert (I-P)^{-1/2}(v_1+v_2)\right\Vert_p.
%\]
As in \cite{asterisque} p. 35 and \cite{ac}, Section 2.4, write
\[
u\varphi=P^{k^2}(u\varphi)+\sum_{l=0}^{k^2-1} P^l(I-P)(u\varphi),
\]
so that
\begin{equation} \label{sumd}
\nabla (u\varphi)\leq \nabla 
\left(P^{k^2}(u\varphi)\right)+\sum_{l=0}^{k^2-1} \nabla \left(P^l(I-P)(u\varphi)\right).
\end{equation}
To treat the first term in the right-hand side of (\ref{sumd}), fix $\rho\in (p,q)$ and notice 
that, since $(R_{\rho})$ holds by assumption, it follows that $k\left\vert 
\nabla P^{k^2}\right\vert$ is $L^{\rho}(\Gamma)$-bounded. Then, arguing as in Lemma \ref{gradball}, one obtains that
\begin{equation} \label{gradballmodif}
\left(\frac 1{V(B)} \sum_{x\in B} \left\vert \nabla P^{l}f(x)\right\vert^pm(x)\right)^{1/p}\leq \frac {Ce^{-\frac{c4^jk^2}l}}{\sqrt{l}} \left(\frac 1{V(2^jB)} \sum_{x\in C_j(B)} \left\vert f(x)\right\vert^2m(x)\right)^{1/2}
\end{equation}
for all $j\geq 1$, all $l\in \left\{1,...,k^2\right\}$ and all function $f$ supported in $C_j(B)$. It follows at once from (\ref{gradballmodif}) applied with $f=u\varphi$, the fact that $u$ has zero integral on $16B$ and the Poincar\'e inequality $(P_2)$ that
\begin{equation} \label{one}
\begin{array}{lll}
\displaystyle \left(\frac 1{V(B)} \sum_{x\in B} \left\vert \nabla P^{k^2}(u\varphi)(x)\right\vert^pm(x)\right)^{1/p}& \leq & \displaystyle \frac Ck\left(\frac 1{V(4B)} \sum_{x\in 4B} \left\vert u(x)\right\vert^2m(x)\right)^{1/2}\\
& \leq & \displaystyle C\left(\frac 1{V(16B)} \sum_{x\in 16B} \left\vert \nabla u(x)\right\vert^2m(x)\right)^{1/2}.
\end{array}
\end{equation}
Let us now turn to the second term in (\ref{sumd}). A calculation shows that, for all $x\in \Gamma$,
\begin{equation} \label{deltaformula}
\begin{array}{lll}
\displaystyle(I-P)(u\varphi)(x) &= & \displaystyle \sum_{y\in \Gamma} p(x,y)((u\varphi)(x)-(u\varphi)(y))\\
& = & \displaystyle \sum_{y\in \Gamma} p(x,y)u(x)(\varphi(x)-\varphi(y))+\sum_{y\in \Gamma} p(x,y) (u(x)-u(y))(\varphi(y)-\varphi(x))\\
& + & \displaystyle \sum_{y\in \Gamma} p(x,y)(u(x)-u(y))\varphi(x)\\
& := & \displaystyle v_1(x)+v_2(x)+v_3(x).
\end{array}
\end{equation}
For all $x\in \Gamma$, $v_3(x)=0$ since $u$ is harmonic in $32B$. Because of the support condition on $\varphi$, one may apply (\ref{gradballmodif}) to $v_2$, and since $\left\Vert \nabla\varphi\right\Vert_{\infty}\leq C/k$, one obtains
\begin{equation} \label{two}
\left(\frac 1{V(B)} \sum_{x\in B} \left\vert \nabla P^l v_2(x)\right\vert^pm(x)\right)^{1/p}\leq \frac C{k\sqrt{l}}\left(\frac 1{V(4B)} \sum_{x\in 4B} \left\vert \nabla u(x)\right\vert^2m(x)\right)^{1/2}
\end{equation}
for all $1\leq l\leq k^2-1$. \par
\noindent For $v_1$, write
\[
\begin{array}{lll}
\displaystyle 2v_1(x) & = & \displaystyle \sum_{y} p(x,y)(u(x)+u(y))(\varphi(x)-\varphi(y))+\sum_{y}
p(x,y)(u(x)-u(y))(\varphi(x)-\varphi(y))\\
& = & \delta F(x)-v_2(x), 
\end{array}
\]
where, for all $(x,y)\in E$,
\[
F(x,y)=(u(x)+u(y))(\varphi(x)-\varphi(y))
\]
is antisymmetric, belongs to $L^2(E)$ and is supported in $B((x_0,x_0),4k)\setminus B((x_0,x_0),2k)$. It is therefore enough to show that, for all $1\leq l\leq k^2-1$,
\begin{equation} \label{doubleoffdiag} 
\left(\sum_{x\in B} \left\vert \nabla P^l \delta F(x)\right\vert^pm(x)\right)^{\frac 1p} \leq \frac {Ce^{-c\frac{k^2}l}}{l} \left(\frac 1{V(B(x_0,x_0),4k)}\sum_{(x,y)\in B((x_0,x_0),4k)\setminus B((x_0,x_0),2k)} \left\vert F(x,y)\right\vert^2 \mu_{xy}\right)^{\frac 12}. 
\end{equation}
To prove this inequality, if $l=2m$, write $\nabla P^l \delta F=\nabla P^m P^m\delta F$. We establish (\ref{doubleoffdiag}) by arguments similar to the proof of Lemma \ref{gradball}, combining (\ref{gradl2off}) and an inequality analogous to (\ref{gradl2off} and derived by duality (see the proof of (2.6) in \cite{ac}). We finally obtain
\begin{equation} \label{three}
\left(\frac 1{V(B)} \sum_{x\in B} \left\vert \nabla P^l v_1(x)\right\vert^pm(x)\right)^{1/p}\leq \frac C{k\sqrt{l}} \left(\frac 1{V(4B)} \sum_{x\in 4B} \left\vert \nabla u(x)\right\vert^2m(x)\right)^{1/2}
\end{equation}
for all $1\leq l\leq k^2-1$. Summing up (\ref{two}) and (\ref{three}) for $1\leq l\leq k^2-1$,  we obtain
\begin{equation} \label{sumup}
\left(\frac 1{V(B)} \sum_{x\in B} \nabla \left(\sum_{l=1}^{k^2-1} \left(P^l(I-P)(u\varphi)\right)\right)(x)^pm(x)\right)^{1/p}\leq C\left(\frac 1{V(16B)} \sum_{x\in 16B} \left\vert \nabla u(x)\right\vert^2m(x)\right)^{1/2}.
\end{equation}
What remains to be treated in (\ref{sumd}) is the term $\nabla (I-P)(u\varphi)$. By (\ref{deltaformula}),
\begin{equation} \label{remains} 
\frac {1}{V(B)^{1/p}}\left\Vert \nabla (I-P)(u\varphi)\right\Vert_{L^p(B)}\leq \frac 1{V(B)^{1/p}}\left\Vert \nabla v_1\right\Vert_{L^p(B)} + \frac 1{V(B)^{1/p}}\left\Vert \nabla v_2\right\Vert_{L^p(B)}.
\end{equation}
Let us first deal with $v_1$. By $(R_p)$,
\[
\left\Vert \nabla v_1\right\Vert_{L^p(\Gamma)}\leq C\left\Vert (I-P)^{1/2}v_1\right\Vert_{L^p(\Gamma)}\leq C\left\Vert v_1\right\Vert_{L^p(\Gamma)},
\]
where the last inequality follows from the $L^p$-boundedness of $(I-P)^{1/2}$ (see \cite{cs}, p. 423 and also \cite{cp7}). But $v_1$ is supported in $4B$ and, for all $x\in 4B$,
\[
\left\vert v_1(x)\right\vert\leq \frac Ck\left\vert u(x)\right\vert.
\]
As a consequence,
\[
\left\Vert v_1\right\Vert_{L^p(\Gamma)}\leq \frac Ck \left\Vert u\right\Vert_{L^p(4B)}\leq \frac Ck \left\Vert u\psi\right\Vert_{L^p(8B)},
\]
where $\psi$ is a nonnegative function equal to $1$ on $4B$, supported in $8B$ and satisfying $\left\Vert \nabla\psi\right\Vert_{\infty}\leq \frac Ck$. Now, (\ref{sob}) shows that, if $q_0=\frac 2{1-\nu}$ and $p\in \left(2,q_0\right)$,
\[
\frac 1{V(B)^{1/p}} \left\Vert u\psi\right\Vert_{L^p(8B)}\leq \frac 1{V(8B)^{1/q_0}} \left\Vert u\psi\right\Vert_{L^{q_0}(8B)}\leq \frac C{V(8B)^{1/q_0}}kV(8B)^{-\nu/2}\left\Vert \nabla (u\psi)\right\Vert_{L^2(8B)}.
\]
Using now the fact that $\frac{\nu}2=\frac 12-\frac 1{q_0}$, we finally conclude
\begin{equation} \label{four}
\frac 1{V(B)^{1/p}} \left\Vert v_1\right\Vert_{L^p(\Gamma)}\leq \frac C{V(8B)^{1/2}} \left\Vert \nabla (u\psi)\right\Vert_{L^2(8B)}\leq \frac C{V(16B)^{1/2}}\left\Vert \nabla u\right\Vert_{L^2(16B)},
\end{equation}
where the last inequality is due $(P_2)$. All these computations yield
\begin{equation} \label{fourbis}
\frac 1{V(B)^{1/p}} \left\Vert \nabla v_1\right\Vert_{L^p(\Gamma)}\leq \frac C{V(16B)^{1/2}}\left\Vert \nabla u\right\Vert_{L^2(16B)}
\end{equation}
\noindent We argue similarly for $v_2$. We just have to notice that, for all $x\in 4B$,
\[
\left\vert v_2(x)\right\vert^p\leq \frac C{k^{p}} \sum_{y\sim x} \left(\left\vert u(y)\right\vert^p+\left\vert u(x)\right\vert^p\right),
\]
hence
\[
\sum_{x\in 4B} \left\vert v_2(x)\right\vert^pm(x)\leq \frac C{k^p} \sum_{x\in 4B} \sum_{y\sim x} \left\vert u(y)\right\vert^p m(x)+ \frac C{k^p} \sum_{x\in 4B} \sum_{y\sim x} \left\vert u(x)\right\vert^p m(x).
\]
Since $m(x)\leq Cm(y)$ whenever $x\sim y$ (this is a straightforward consequence of (\ref{D}) and was noticed in \cite{cougri}, Section 4.2) and $\sharp\left\{y\in \Gamma;\ y\sim x\right\}\leq N$, we finally obtain that
\[
\sum_{x\in 4B} \left\vert v_2(x)\right\vert^pm(x)\leq \frac C{k^p} \sum_{x\in 8B} \left\vert u(x)\right\vert^pm(x),
\]
and we conclude as for $v_1$ that
\begin{equation} \label{five}
\frac 1{V(B)^{1/p}} \left\Vert \nabla v_2\right\Vert_{L^p(\Gamma)}\leq \frac C{V(16B)^{1/2}}\left\Vert \nabla u\right\Vert_{L^2(16B)},
\end{equation}
Summing up (\ref{one}), (\ref{sumup}, (\ref{fourbis}) and (\ref{five}), we obtain that $(RH_p)$ holds.\hfill\fin\par
\noindent As far as Proposition \ref{selfimprove} is concerned, its proof is entirely similar to the one of Proposition 2.2 in \cite{ac} and will therefore be skipped. Let us just mention that it relies on an elliptic Caccioppoli inequality (analogous to the Euclidean version, see \cite{gia}), Proposition \ref{poincareself} and Gehring's self-improvement of reverse H\"older inequalities (\cite{gehring}).\hfill\fin
\section*{Appendix} \label{app}\setcounter{equation}{0}
\addcontentsline{toc}{section}{Appendix}
We prove Lemma \ref{estim}. For all $l\geq 0$, $\displaystyle a_l=\frac{(2l)!}{4^l(l!)^2}$, and, as already used in Section \ref{reverseriesz}, the Stirling formula shows $a_l\sim\frac 1{\sqrt{\pi l}}$. Therefore, there exists $C>0$ such that, for all $l\geq 1$,
\[
0<a_l\leq \frac C{\sqrt{l}}.
\]
Assume first that $mk^2<l<(m+1)k^2$ for some integer $0\leq m\leq n$. For each integer $j\geq 0$ such that $jk^2\leq l$, one has $l-jk^2>0$ and $j\leq m$, so that $\left\vert a_{l-jk^2}\right\vert\leq \frac C{\sqrt{l-jk^2}}\leq \frac C{\sqrt{l-mk^2}}$. It follows at once that
\[
\left\vert d_l\right\vert\leq \frac C{\sqrt{l-mk^2}}
\]
for some $C>0$ only depending on $n$. \par
Assume now that $l=(m+1)k^2$ for some $0\leq m\leq n$. For each $j\geq 0$ such that $jk^2\leq l$ and $l-jk^2>0$, one has $j\leq m$ again, so that $\left\vert a_{l-jk^2}\right\vert\leq \frac C{\sqrt{l-mk^2}}=\frac Ck\leq C$. Moreover, $a_0=1$. One therefore has
\[
\left\vert d_l\right\vert\leq C+C_n^{m+1}\leq C,
\]
where, again, $C$ only depends on $n$.\par
Finally, assume that $l>(n+1)k^2$. The classical computation of Wallis integrals shows that
\[
a_l=\frac 2{\pi}\int_0^{\frac{\pi}2} \left(\sin t\right)^{2l}dt=\varphi(l)
\]
where, for all $x>0$, $\displaystyle \varphi(x)=\frac 2{\pi}\int_0^{\frac{\pi}2} \left(\sin t\right)^{2x}dt$. We can then invoke (\ref{fellerineq}) and are therefore left with the task of estimating $\varphi^{(n)}$. But, for all $x>0$,
\[
\left\vert \varphi^{(n)}(x)\right\vert =\frac 2{\pi} \left\vert \int_0^{\frac{\pi}2} \left(2\log\sin t\right)^ne^{2x\log\sin t}dt\right\vert \leq \frac 2{\pi} \int_0^{\frac{\pi}2} \left\vert 2\log\sin t\right\vert^ne^{2x\log\sin t}dt:=\frac 2{\pi}I_n(x).
\]
We now argue as in the ``Laplace'' method. For all $\delta\in \left(0,\frac{\pi}2\right)$, one clearly has, for all $x>1$,
\begin{equation} \label{onebis}
\begin{array}{lll}
0\leq I_n(x)& \leq &\displaystyle \int_0^{\frac{\pi}2-\delta} \left\vert 2\log\sin t\right\vert^ne^{2x\log\sin t}dt+\int_{\frac{\pi}2-\delta}^{\frac{\pi}2} \left\vert 2\log\sin t\right\vert^ne^{2x\log\sin t}dt\\
& \leq & \displaystyle \left(\sin\left(\frac{\pi}2-\delta\right)\right)^{2x-2} I_n(1)+J_n(x)=C_{n,\delta}\alpha^{2x-2}+ J_n(x)
\end{array}
\end{equation}
where $C_{n,\delta}>0$ only depends on $n$ and $\delta$, $0<\alpha=\sin\left(\frac{\pi}2-\delta\right)<1$ and $\displaystyle J_n(x):=\int_{\frac{\pi}2-\delta}^{\frac{\pi}2} \left\vert 2\log\sin t\right\vert^ne^{2x\log\sin t}dt$. \par
Observe now that $\displaystyle J_n(x)=\int_0^{\delta} \left\vert 2\log\cos u\right\vert^n e^{2x\log\cos u}du$. Since $\log(\cos u)\sim-\frac{u^2}2$ when $u\rightarrow 0$, we fix $\delta>0$ such that, for all $0<u<\delta$, $-\frac 34 u^2\leq \log(\cos u)\leq -\frac 14 u^2$, which implies
\begin{equation} \label{twobis}
\begin{array}{lll}
J_n(x) & \leq & \displaystyle C \int_0^{\delta} u^{2n} e^{-\frac 12 xu^2}du\\
& \leq & \displaystyle C \left(\frac 1{\sqrt{x}}\right)^{2n+1}\int_0^{+\infty} v^{2n}e^{-v^2}dv\leq Cx^{-n-\frac 12}.
\end{array}
\end{equation}
It follows from (\ref{onebis}) and (\ref{twobis}) that, for all $x>1$,
\[
\left\vert \varphi^{(n)}(x)\right\vert\leq Cx^{-n-\frac 12},
\]
which, joined with (\ref{fellerineq}), yields assertion $(iii)$ in Lemma \ref{estim}, the proof of which is now complete. \hfill\fin


\begin{thebibliography}{AAA}
%\bibitem{alexo} G. Alexopoulos, An application of homogenization theory 
%to harmonic analysis: Harnack inequalities and Riesz transforms on Lie 
%groups of polynomial volume growth, {\it Can. J. Math.}, 44, 4, 
%691-727, 1992.
\bibitem{a} P. Auscher, On $L^p$ estimates for square roots of second order elliptic operators on $\R^n$, {\it Publ. Mat.}, 48, 159-186, 2004.
\bibitem{auscher} P. Auscher, On necessary and sufficient conditions 
for $L^p$ estimates of Riesz transforms associated to elliptic 
operators on $\R^n$ and related estimates, {\it Mem. A. M. S.}, vol. 
186, 871, March 2007.
\bibitem{ac} P. Auscher, T. Coulhon, Riesz transforms on manifolds 
and Poincar\'e inequalities, {\it  Ann. Sc. Norm. Super. Pisa Cl. 
Sci.} (5)  4,  no. 3, 531--555, 2005.
\bibitem{auschcou} P. Auscher, T. Coulhon, Gaussian lower bounds for random 
walks from elliptic regularity, {\it Ann. Inst. H. Poincar\'e Probab. Statist.} 
35, 5, 605--630, 1999.
\bibitem{acdh} P. Auscher, T. Coulhon, X. T. Duong, S. Hofmann, Riesz 
transforms 
on manifolds and heat kernel regularity, {\it Ann. Sci. Ecole Norm. 
Sup.}, 37, 6, 911-957, 2004.
\bibitem{am} P. Auscher, J.-M. Martell, Weighted norm inequalities, 
off-diagonal estimates and elliptic operators: Part IV, preprint, 
http://arxiv.org/format/math.DG/0603643, 2006.
\bibitem{asterisque} P. Auscher, P. Tchamitchian,{\it Square root
problem for divergence operators and related topics}, Ast\'erisque,
{\bf 249}, Soc. Math. France (1998).
\bibitem{badr1} N. Badr, Real interpolation of Sobolev spaces, preprint, http://arxiv.org/abs/0705.2216v2.
\bibitem{bakry} D. Bakry, Etude des transformations de Riesz dans les 
vari\'et\'es riemanniennes \`a courbure de Ricci minor\'ee, in {\it 
S\'eminaire de probabilit\'es XXI}, Springer L. N. n¡ 1247, 137-172, 
1987.
\bibitem{bcls} D. Bakry, T. Coulhon, M. Ledoux, L. Saloff-Coste, Sobolev inequalities in disguise, {\it Indiana Univ. Math. J.} 124, 4, 1033-1074, 1995. 
\bibitem{bs} C. Bennett, R. Sharpley, {\it Interpolation of operators}, Academic Press, 1988.
\bibitem{blunck} S. Blunck, Perturbation of analytic operators and 
temporal regularity of discrete heat kernels, {\it Coll. Math.} 86, 
189-201, 2000.
\bibitem{bk} S. Blunck, P. Kunstmann, Calder\'on-Zygmund theory for non-integral operators and the $H^{\infty}$ functional calculus, 
{\it Rev. Mat. Iberoamericana} 19, 3, 919-942, 2003.
\bibitem{cks} E. Carlen, S. Kusuoka, D. Stroock, Upper bounds for symmetric Markov transition functions, {\it Ann. Inst. H. Poincar\'e, proba. stat.}, suppl. au no 2, 245-287, 1987.
\bibitem{carron} G. Carron, In\'egalit\'es isop\'erim\'etriques de Faber-Krahn et cons\'equences, in {\it Actes de la table ronde de g\'eom\'etrie diff\'erentielle en l'honneur de Marcel Berger}, Coll. SMF S\'eminaires et Congr\`es, 1, 1994.
\bibitem{christ} M. Christ, Temporal regularity for random walk on 
discrete nilpotent groups, {\it J. Fourier Anal. Appl.}, Kahane 
special issue, 141-151, 1995.
\bibitem{coifw} R. Coifman, G. Weiss, {\it Analyse harmonique 
non-commutative sur certains espaces homog\`enes}, Springer L. N. 242, 
1971.
\bibitem{coifman1} R. Coifman, G. Weiss, Extensions of Hardy
spaces and their use in analysis, {\it Bull. Amer. Math. Soc.} 83, 
569-645, 1977.
\bibitem{cp7} T. Coulhon, Suite d'op\'erateurs \`a puissances born\'ees dans les espaces ayant la propri\'et\'e de Dunford-Pettis, {\it S\'eminaire d'analyse fonctionnelle 84/85}, Publications math\'ematiques de l'Universit\'e Paris VII, 26, 1986.
\bibitem{coulhonlp} T. Coulhon, Espaces de Lipschitz et in\'egalit\'es de Poincar\'e, {\it J. Funct. Anal.} 136, 1, 81-113, 1996. 
\bibitem{coulhonsurvey} T. Coulhon, Random walks and geometry on infinite graphs,  in {\it Lecture notes on analysis on metric spaces, Trento, C.I.R.M., 1999}, Luigi Ambrosio, Francesco Serra Cassano, ed., Scuola Normale Superiore di Pisa, 5-30, 2000.
\bibitem{cd} T. Coulhon, X. T. Duong, Riesz transforms for $1\leq 
p\leq 2$, {\it Trans. Amer. Math. Soc.}  351, no. 3, 1151-1169, 1999.
\bibitem{cougri} T. Coulhon, A. Grigor'yan, Random walks on graphs with
regular volume growth, {\it Geom. and Funct. Anal.} 8, 656-701, 1998.
\bibitem{cl} T. Coulhon, M. Ledoux, Isop\'erim\'etrie, d\'ecroissance du noyau de la chaleur et transformations de Riesz: un contre-exemple, {\it Ark. f\"or Mat.} 32, 63-77, 1994.
\bibitem{cs} T. Coulhon, L. Saloff-Coste, Puissances d'un op\'erateur 
r\'egularisant,  {\it Ann. Inst. H. Poincar\'e Probab. Statist.}  26,  
no. 3, 419--436, 1990.
%\bibitem{davies1} E. B. Davies, Heat kernel bounds, conservation of 
%probability and the Feller property, {\it J. Analyse Math.} 58, 
%99-119, 1992.
%\bibitem{davies2} E. B. Davies, Uniformly elliptic operators with 
%measurable 
%coefficients, {\it J. Funct. Anal.}  132,  no. 1, 141--169, 1995.
\bibitem{delmoelliptic} T. Delmotte, In\'egalit\'e de Harnack elliptique sur les graphes, {\it Coll. Math.} 72, 1, 19-37, 1997.
\bibitem{delmo} T. Delmotte, Parabolic Harnack inequality, {\it Rev. Mat. Iberoamericana}, 15, 1, 181-232, 1999. 
\bibitem{dungey} N. Dungey, A note on time regularity for discrete time heat 
kernels, {\it Semigroup Forum}, 72, 404-410, 2006. 
\bibitem{dungeylp} N. Dungey, A Littewood-Paley-Stein estimate on graphs and groups, preprint.
\bibitem{feller} W. Feller, {\it An introduction to probability theory and its applications}, volume I, Wiley, 1968.
\bibitem{gaffney} M. P. Gaffney, The conservation property of the 
heat 
equation on Riemannian manifolds, {\it Comm. Pure Appl. Math.} 12, 
1-11, 1959.
\bibitem{gehring} F. W. Gehring, The $L^p$ integrability of the partial derivative of a quasi-conformal mapping, {\it Acta Math.} {\bf 130}, 265-277, 1973. 
\bibitem{gia} M. Giaquinta, {\it Multiple integrals in the calculus of variations and non-linear elliptic equations}, volume 105 of Annals of Math. Studies, Princeton Univ. Press, 1983.
\bibitem{gt} V. Gol'dshtein and M. Troyanov, Axiomatic theory of Sobolev Spaces, {\it Expo. Math.}, 19, 289-336, 2001.
\bibitem{grigo} A. Grigor'yan, Heat kernel upper bounds on a complete non-compact manifold, {\it Rev. Mat. Iberoamericana}  10, 2, 395-452, 1994.
\bibitem{hk} P. Hajlasz, P. Koskela, Sobolev met Poincar\'e,  Mem. Amer. 
Math. Soc.  145,  no. 688, 2000.
\bibitem{hebsal} W. Hebisch, L. Saloff-Coste, Gaussian estimates for Markov 
chains and random walks on groups, {\it Ann. Probab.} 21, 2, 
673--709, 1993.
\bibitem{hm} S. Hofmann, J.-M. Martell, $L^p$ bounds for Riesz transforms and square roots associated to second order elliptic divergence operators, {\it Publ. Mat.} 47, 497-515, 2003.
\bibitem{kz} S. Keith, X. Zhong, The Poincar\'e inequality is an open 
ended condition, preprint of the University of Jyv\"askyl\"a, 2003, to appear in {\it Ann. Math.}.
\bibitem{martellthesis} J. M. Martell, Desigualdades con pesos en el 
An\'alisis 
de Fourier: de los espacios de tipo homog\'eneo a las medidas no 
doblantes, 
Ph. D., Universidad Aut\'onoma de Madrid, 2001.
\bibitem{meyers} N. G. Meyers, An $L^p$ estimate for the gradient of solutions of second order elliptic divergence equations, {\it Ann. Scuola Norm. Sup. Pisa} {\bf 3}, 17, 189-206, 1963.
\bibitem{ostr} M. I. Ostrovskii, Sobolev spaces on graphs, {\it Quaest. Math.}, 28, 4, 501-523, 2005.
%\bibitem{theserey} T. Rey, Estimations de De Giorgi Nash et approximations, Th\`ese de doctorat de l'Universit\'e Paul C\'ezanne, 2004.
\bibitem{r} E. Russ, Riesz transforms on graphs for $1 \leq p \leq 2$, {\it Math.Scand.}, 87, 133-160, 2000.
\bibitem{pota} E. Russ, $H^1-L^1$ boundedness of Riesz transforms on 
Riemannian manifolds and on graphs, {\it Pot. Anal.}, 14, 301-330, 
2001.
\bibitem{shen} Z. Shen, Bounds of Riesz transforms on $L^p$ spaces for second order elliptic operators, {\it Ann. Inst. Fourier} 55, 1, 173-197, 2005.
\bibitem{topics} E. M. Stein, Topics in harmonic analysis related to the Littlewood-Paley theory, Princeton Univ. Press, 1970.
%\bibitem{stein} E. M. Stein, {\it Harmonic analysis: real-variable methods, 
%orthogonality and oscillatory integrals}, Princeton Univ. Press, 1993.
\bibitem{stri} R. Strichartz, Analysis of the Laplacian on the 
complete Riemannian manifold, {\it J. Funct. Anal.}, 52, 1, 48-79, 1983.
%\bibitem{varo} N. Varopoulos, Isoperimetric inequalities and Markov chains, {\it J. Funct. Anal.} 63, 215-239, 1985.
\end{thebibliography}
\end{document}